\documentclass[11pt,leqno]{amsart}
\usepackage[mathcal]{eucal}
\usepackage[utf8x]{inputenc}
\usepackage[T1]{fontenc}
\usepackage[english]{babel}

\usepackage{amsmath,amssymb,amsthm}
\usepackage[dvipsnames]{xcolor}
\usepackage[colorlinks=true, linkcolor=MidnightBlue, citecolor=OliveGreen, pagebackref]{hyperref}
\usepackage[lite,initials,alphabetic,backrefs]{amsrefs}
\usepackage[capitalise,nameinlink]{cleveref}
\usepackage{enumerate}

\renewcommand{\PrintReviews}[1]{}

\PrerenderUnicode{ç}

\numberwithin{equation}{section}

\newcounter{daggerfootnote}
\newcommand*{\daggerfootnote}[1]{%
    \setcounter{daggerfootnote}{\value{footnote}}%
    \renewcommand*{\thefootnote}{\fnsymbol{footnote}}%
    \footnote[2]{#1}%
    \setcounter{footnote}{\value{daggerfootnote}}%
    \renewcommand*{\thefootnote}{\arabic{footnote}}%
    }

\newtheorem{theorem}{Theorem}[section]
\newtheorem{lemma}[theorem]{Lemma}
\crefname{lemma}{Lemma}{Lemmata}
\newtheorem{proposition}[theorem]{Proposition}
\newtheorem{corollary}[theorem]{Corollary}

\newtheorem{problem}{Problem}

\theoremstyle{definition}

\newtheorem{example}[theorem]{Example}
\newcommand{\Z}{\mathbb{Z}}
\newcommand{\Q}{\mathbb{Q}}
\newcommand{\N}{\mathbb{N}}
\newcommand{\nt}{\trianglelefteq}

\DeclareMathOperator{\Zent}{Z}
\DeclareMathOperator{\cc}{C^*\!\!}

\title[Free polynilpotent groups and the Magnus property]{Free polynilpotent groups and \\ the Magnus property}

\author[B. Klopsch]{Benjamin Klopsch}
\address{Benjamin Klopsch:
	Mathematisches Institut, Heinrich-Heine-Universit\"at, 40225
	D\"usseldorf, Germany}
\email{klopsch@math.uni-duesseldorf.de}

\author[L. Mendonça]{Luis Mendonça}
\address{Luis Mendonça:
  Instituto de Ci\^ encias Exatas, Universidade Federal de Minas Gerais, 31270-901
  Belo Horizonte, Brazil}
\email{demendoncaluisaugusto@gmail.com}

\author[J. M. Petschick]{Jan Moritz Petschick}
\address{Jan Moritz Petschick:
	Mathematisches Institut, Heinrich-Heine-Universit\"at, 40225
	D\"usseldorf, Germany}
\email{jan.petschick@hhu.de}

\thanks{The research was partially conducted in the framework of the
  DFG-funded research training group “GRK 2240: Algebro-Geometric
  Methods in Algebra, Arithmetic and Topology”. Furthermore, the
  research was partially funded by the Deutsche Forschungsgemeinschaft
  (DFG, German Research Foundation) — 380258175. The second author was
  supported by the S\~ao Paulo Research Foundation (FAPESP) grant
  21/13052-5.}


 
\begin{document}

\begin{abstract}
  Motivated by a classic result for free groups, one says that a group
  $G$ has the Magnus property if the following holds: whenever two
  elements generate the same normal subgroup of~$G$, they are
  conjugate or inverse-conjugate in~$G$.
  
  It is a natural problem to find out which relatively free groups
  display the Magnus property.  We prove that a free polynilpotent
  group of any given class row has the Magnus property if and only if
  it is nilpotent of class at most~$2$.  For this purpose we explore
  the Magnus property more generally in soluble groups, and we produce
  new techniques, both for establishing and for disproving the
  property.  We also prove that a free centre-by-(polynilpotent of
  given class row) group has the Magnus property if and only if it is
  nilpotent of class at most~$2$.

  On the way, we display $2$-generated nilpotent groups (with
  non-trivial torsion) of any prescribed nilpotency class with the
  Magnus property.  Similar examples of finitely generated,
  torsion-free nilpotent groups are hard to come by, but we construct
  a $4$-generated, torsion-free, class-$3$ nilpotent group of Hirsch
  length $9$ with the Magnus property.  Furthermore, using a weak
  variant of the Magnus property and an ultraproduct construction, we
  establish the existence of metabelian, torsion-free, nilpotent
  groups of any prescribed nilpotency class with the Magnus property.
\end{abstract}

\maketitle

%%%%%

\section{Introduction}\label{sec:introduction}

A group $G$ has the \emph{Magnus property} if the following holds:
whenever $g,h \in G$ generate the same normal subgroup
$\langle g \rangle^G = \langle h \rangle^G$, the element $g$ is
already conjugate in~$G$ to $h$ or to~$h^{-1}$.  Magnus~\cite{Ma30}
established this property for free groups, using his
``Freiheitssatz''.  The Magnus property is a first-order property, in
the sense of model theory; consequently all groups with the same
elementary theory as free groups have the Magnus property.  During the
last two decades, the Magnus property has been explored and
established for various classes of groups using different techniques,
e.g., for fundamental groups of closed surfaces~\cite{BoSv08}, direct
products of free groups~\cites{KlKu16,Fe21}, and certain amalgamated
products \cite{Fe19}.

Groups with the Magnus property are typically torsion-free and `big',
for instance, in the sense that they do not satisfy any non-trivial
law, viz.\ any non-trivial identical relation.  Even so free abelian
groups possess the Magnus property, for obvious reasons, and certain
crystallographic groups with the Magnus property were manufactured
in~\cite{KlKu16}.  In conjunction with Magnus' classic result, this
prompts a natural question for relatively free groups, viz.\
$\mathcal{V}$-free groups for any given variety $\mathcal{V}$ of
groups.

\begin{problem} \label{prob:A}
  Let $\mathcal{V}$ be a variety of groups, viz.\ the class of all
  groups satisfying each one of a given set of laws.  Which
  $\mathcal{V}$-free groups have the Magnus property?
\end{problem}

By basic considerations, it is enough to settle the question for
relatively free groups of finite rank, viz.\ on finitely many free
generators, because for each variety of groups~$\mathcal{V}$ there is
a precise cut-off point
$\delta_\mathcal{V} \in \N_0 \cup \{ \infty \}$ for the ranks of
$\mathcal{V}$-free groups with the Magnus property; see
\cref{cor:f-g-enough}.  As mentioned above, for the variety
$\mathcal{U}$ of all groups and for the variety $\mathcal{A}$ of
abelian groups, we know that every free group and every
$\mathcal{A}$-free group has the Magnus property: hence
$\delta_\mathcal{U} = \delta_\mathcal{A} = \infty$.  Likewise, it is
easy to see that the variety $\mathcal{A}_m$ of abelian groups of
exponent $m$ (that is, $m$ or dividing $m$) has
$\delta_{\mathcal{A}_m} = \infty$ if $m \in \{1,2,3,4,6\}$, and
$\delta_{\mathcal{A}_m} = 0$ otherwise.

Perhaps it is natural to concentrate first on varieties of exponent
zero, viz.\ varieties $\mathcal{V}$ such that $x^m$ is not a universal
law in $\mathcal{V}$-groups, for any~$m \in \N$.  Prominent examples
of this kind are the varieties $\mathcal{N}_\mathbf{c}$ of all
\emph{polynilpotent groups} of class row $\mathbf{c}$, for any given
length $l \in \N$ and class tuple
$\mathbf{c} = (c_1,\ldots,c_l) \in \N^l$.  We recall that a group $G$
belongs to $\mathcal{N}_\mathbf{c}$ if the term
$\gamma_{(c_1+1,\ldots,c_l+1)}(G)$ of its \emph{iterated lower central
  series} vanishes; here $\gamma_{(1)}(G) = \gamma_1(G) = G$ and
inductively we set
\[
  \gamma_{(c_1+1,\ldots,c_l+1)}(G) = \gamma_{c_l+1} \big(
  \gamma_{(c_1+1,\ldots,c_{l-1}+1)}(G) \big), \qquad \text{for $l>1$},
\]
where $\gamma_{(c_1+1)}(G) = \gamma_{c_1+1}(G) = [\gamma_{c_1}(G),G]$
is the $(c_1+1)$\textsuperscript{st} term of the ordinary lower
central series of~$G$.  For instance, for $l=1$ and $\mathbf{c} = (c)$
the variety $\mathcal{N}_\mathbf{c}$ consists of all nilpotent groups
of class at most~$c$; for $l \in \N$ and
$\mathbf{c} = (1,\ldots,1) \in \N^l$, the variety
$\mathcal{N}_\mathbf{c}$ consists of all soluble groups of derived
length at most~$l$.  For free polynilpotent groups, we resolve
\cref{prob:A} completely.

\begin{theorem}\label{thm:free-polynilpotent}
  Let $G$ be an $\mathcal{N}_\mathbf{c}$-free group of rank~$d$, i.e.,
  a free polynilpotent group of class row~$\mathbf{c}$ that is freely
  generated by $d$ elements, where $d,l \in \N$ and
  $\mathbf{c} \in \N^l$.

  Then $G$ has the Magnus property if and only if $G$ is nilpotent of
  class at most~$2$; equivalently, if and only if $d=1$ or
  $\mathbf{c} \in \{ (1), (2) \}$.
\end{theorem}

The proof uses the notion of basic witness pairs for not having the
Magnus property; which are defined in~\cref{lem:basic-witness-pair}.
The starting point of the proof is that the restricted wreath product
$C_\infty \wr C_\infty$ admits such witness pairs; see
\cref{prop:wreath-product}. Almost as a by-product, we obtain the
following similar result, for further varieties.

\begin{theorem} \label{thm:free-centre-by-polynilpotent}
  Let $G$ be a free centre-by-$\mathcal{N}_\mathbf{c}$ group of
  rank~$d$, where $d,l \in \N$ and $\mathbf{c} \in \N^l$.  Then $G$
  has the Magnus property if and only if $G$ is nilpotent of class at
  most~$2$; equivalently, if and only if $d=1$ or
  $\mathbf{c} = (1)$.
\end{theorem}

Since the Magnus property is a first-order property, one may wonder
about groups with the same elementary theory as free polynilpotent
groups.  Groups that are elementarily equivalent to free nilpotent
groups were considered in~\cites{MySo09,MySo11}.  We remark that
$\mathcal{N}_\mathbf{c}$-free groups are torsion-free, while free
centre-by-$\mathcal{N}_\mathbf{c}$ groups may involve central torsion
of exponent~$2$; compare with~\cites{Ku82,St89}.

\medskip

In order to prove
\cref{thm:free-polynilpotent,thm:free-centre-by-polynilpotent} we
explore the Magnus property in more general groups, and we produce new
techniques, both for establishing and for disproving the property.  In
\cref{prop:quot-crit} we provide a useful sufficient criterion under
which the Magnus property passes to factor groups; for instance, if
$G$ has the Magnus property and $N \nt G$ is finite then $G/N$
inherits the Magnus property.  In \cref{prop:class-2-nilp->mp} we see
that every torsion-free, class-$2$ nilpotent group has the Magnus
property.  \cref{exa:fg-nilp-MP-any-class} provides an explicit family
of finitely generated, nilpotent groups (with non-trivial $3$-torsion)
of any prescribed nilpotency class that possess the Magnus property.
In contrast it appears much harder to capture finitely generated,
nilpotent groups of prescribed nilpotency class $c \ge 3$ with the
Magnus property that are torsion-free.  In
\cref{exa:4-gen-class-3-MP-example} we construct explicitly a
$4$-generated, torsion-free, class-$3$ nilpotent group of Hirsch
length~$9$ with the Magnus property.  This can be seen as a very first
step toward tackling the following problem which suggests itself.

\begin{problem}
  Are there finitely generated, torsion-free, nilpotent groups with the
  Magnus property of any prescribed nilpotency class?  Characterise or
  even classify finitely generated, torsion-free, nilpotent groups with
  the Magnus property.
\end{problem}

Currently, we seem to be far from solving this problem, but we can
establish a related and somewhat surprising result, using a weak
variant of the Magnus property and an ultraproduct construction.

\begin{theorem}\label{thm:nilpotent-existence}
  For every $c \in \N$ there exists a countable, metabelian,
  torsion-free, nilpotent group with the Magnus property that has
  nilpotency class precisely~$c$.
\end{theorem}

In this context we remark that, by a classical embedding theorem of
Higman, Neumann and Neumann~\cite{HiNeNe49}, every countable
torsion-free group $G$ can be embedded into a countable torsion-free
group $\mathcal{G}$ with only two conjugacy classes; such a group
$\mathcal{G}$ has the Magnus property.  Using small cancellation
techniques, Osin~\cite{Os10} showed that one can even arrange for
$\mathcal{G}$ to be finitely generated.  However, the structure of
such groups $\mathcal{G}$, which arise as inductive limits, can be
very different from the one of the input group~$G$, and they are far
from the groups we produce for \cref{thm:nilpotent-existence}.

\medskip  

\noindent \emph{Notation.}  Let $G$ be a group.  For $x,y \in G$ we
write $x^y = y^{-1}xy$ and $[x,y] = x^{-1} x^y$.  Throughout, we use
left-normed commutators; for instance, we write $[x,y,z] = [[x,y],z]$.
A similar convention applies to iterated Lie commutators in associated
Lie rings.  For $X \subseteq G$ we denote by
$\langle X \rangle^G = \langle x^g \mid x \in X,\, g \in G \rangle$
the normal closure of $X$ in~$G$, viz.\ the smallest normal subgroup
containing~$X$.  For a singleton $X = \{x\}$ we use the shorter
notation $\langle x \rangle^G$.

We denote by $\Zent(G)$ the centre of~$G$, and we write $\Zent_i(G)$,
$i \in \N_0$, for the terms of the upper central series of~$G$.  The
iterated lower central series and, in particular, the lower central
series $\gamma_i(G)$, $i \in \N$, were already discussed above.  

Suppose that $\mathcal{V}$ is a non-trivial variety of groups.  The
\emph{rank} of a $\mathcal{V}$-free group $G$ is the cardinality of a
$\mathcal{V}$-free generating set for $G$.  We use the term sparingly
and no confusion with other common notions of rank, such as Pr\"ufer
rank should arise.

For $m \in \N \cup \{ \infty \}$ we write $C_m$ to denote a cyclic
group of order~$m$.

%%%
\section{Preliminaries and auxiliary results}

We recall that the \emph{Magnus property} is a first-order property in
the sense of model theory; indeed, sometimes it is useful to rephrase
it for a group $G$ as follows:
\begin{equation} \label{equ:magnus-prop} \tag{$\mathsf{MP}$}
  \begin{split} 
    \forall k,l \in \N_0 \quad \forall g,h \in G  \qquad & \forall
    m_1, \ldots, m_k \in \{1,-1\} \quad \forall
    v_1, \ldots, v_k \in G \\
    \qquad & \forall n_1, \ldots, n_l \in \{1,-1\} \quad \forall
    w_1, \ldots, w_l \in G : \\
     \Big( h = \prod\nolimits_{i=1}^k (g^{m_i})^{v_i} \; \land \; g =
    \prod\nolimits_{j=1}^l & (h^{n_j})^{w_j} \Big)
    \implies \Big( \exists v \in G:\quad g^v = h \; \lor \;
    g^v = h^{-1} \Big),
  \end{split}
\end{equation}
where the quantifier over the integers $k,l$ can be eliminated by
passing to a countable collection of sentences in the first-order
language of groups.  For short we say that $G$ is an
\emph{$\mathsf{MP}$-group} if $G$ has the Magnus
property.\protect\daggerfootnote{The terms ``$M$-group'' and ``Magnus
  group'' are unfortunately already in use with other meanings.  For
  lack of better alternatives, we have settled for
  ``$\mathsf{MP}$-group''.}

We recall that, if $\mathcal{P}$ is any property of groups, then a
group $G$ is \emph{locally a $\mathcal{P}$-group} if each finite
subset of $G$ is contained in a $\mathcal{P}$-subgroup of~$G$.  If
$\mathcal{P}$ is inherited by subgroups, this is equivalent to the
requirement that each finitely generated subgroup of $G$
has~$\mathcal{P}$.  The proof of the following lemma is routine,
using~\eqref{equ:magnus-prop}.

\begin{lemma} \label{lem:locally-finite} Every locally
  $\mathsf{MP}$-group is an $\mathsf{MP}$-group.
\end{lemma}

Of course, the Magnus property does not generally pass from a group to
its subgroups or quotients.  Nevertheless there are interesting
situations, where this happens.  We record a simple, but useful
observation.

\begin{lemma} \label{lem:retract} Let $H \le G$ be a retract of a
  group~$G$, that is $G = H \ltimes N$ for a normal complement
  $N \trianglelefteq G$.  If $G$ is an $\mathsf{MP}$-group
  then so is $H$.
\end{lemma}

For each variety of groups $\mathcal{V}$ we set
\[
  \delta_\mathcal{V} = \sup \{ d \in \N_0 \mid
  \text{$\mathcal{V}$-free groups of rank $d$ are
    $\mathsf{MP}$-groups} \} \in \N_0 \cup \{ \infty \}.
\]
\cref{lem:locally-finite,lem:retract} already have a useful
consequence for relatively free groups.

\begin{corollary} \label{cor:f-g-enough} Let $\mathcal{V}$ be a
  variety of groups.  If $\delta_\mathcal{V} = \infty$, then every
  $\mathcal{V}$-free group is an $\mathsf{MP}$-group.  If
  $\delta_\mathcal{V} < \infty$, then a $\mathcal{V}$-free group is an
  $\mathsf{MP}$-group if and only if it has rank at most
  $\delta_\mathcal{V}$.
\end{corollary}

The next result is less obvious, if not surprising; in particular, it
provides a powerful handle to deal with free nilpotent and, more
generally, free polynilpotent groups.

\begin{proposition}\label{prop:quot-crit}
  Let $G$ be an $\mathsf{MP}$-group, and let $N \nt G$ such that for
  each $g \in G\smallsetminus N$ the $\subseteq$-partially ordered set
  of normal subgroups
  \[
    \Omega_{gN} = \big\{ \langle gz \rangle^G \mid z \in N \big\}
  \]
  satisfies the minimal condition.  Then $G/N$ is an
  $\mathsf{MP}$-group.
\end{proposition}

\begin{proof}
  Let $g, h \in G$ such that their images in $G/N$ have the same
  normal closure, in other words such that
  $\langle g \rangle^G \equiv_N \langle h \rangle^G$.  If
  $g \equiv_N 1$, also $h \equiv_N 1$, and they are conjugate to one
  another modulo~$N$.  Now suppose that $g \not\equiv_N 1$.  Choose
  $k,l \in \N$, $m_1, \dots, m_k, n_1, \dots, n_l \in \Z$ and
  $v_1, \dots, v_k, w_1, \dots, w_l \in G$ such that
  \begin{align*}
    \prod\nolimits_{i = 1}^k \left( g^{m_i} \right)^{v_i} \equiv_N h
    \qquad \text{and} \qquad
    \prod\nolimits_{j = 1}^l \left( h^{n_j} \right)^{w_j} \equiv_N  g.
  \end{align*}
  Since $\Omega_{gN}$ satisfies the minimal condition, we find
  $g_{\rm{min}} \in g N$ such that $\langle g_{\rm{min}} \rangle^G$ is
  $\subseteq$-minimal among all subgroups of the form
  $\langle y \rangle^G$ for $y \in gN$.  Consider
  \begin{align*}
    h_0 = \prod\nolimits_{i = 1}^k \left( g_{\rm{min}}^{\, m_i} \right)^{v_i}
    \equiv_N h \qquad \text{and} \qquad
    g_0 = \prod\nolimits_{j = 1}^l \left( h_0^{\, n_j} \right)^{w_j} \equiv_N g.
  \end{align*}
  These elements satisfy
  $\langle g_0 \rangle^G \subseteq \langle h_0 \rangle^G \subseteq
  \langle g_{\rm{min}} \rangle^G$, hence, by the minimal choice of
  $g_{\rm{min}}$, we conclude that
  $\langle g_0 \rangle^G = \langle h_0 \rangle^G$.  Since $G$ has the
  Magnus property, there exists $v \in G$ such that $g_0^{\, v} = h_0$
  or $g_0^{\, v} = h_0^{\, -1}$, hence
  \[
    g^v \equiv_N g_0^{\, v} \equiv_N h_0 \equiv_N = h \quad \text{or,
      similarly,} \quad g^v \equiv_N h^{-1}. \qedhere
  \]
\end{proof}

We record some immediate consequences, which are quite
remarkable.

\begin{corollary}
  Let $G$ be an $\mathsf{MP}$-group,
  $\Omega_G = \{ \langle g \rangle^G \mid g \in G \}$, and let
  $N \nt G$.
  
  $\mathrm{(1)}$ If $\Omega_G$ satisfies the minimal condition, then $G/N$
  is an $\mathsf{MP}$-group.

  $\mathrm{(2)}$ If $N$ is finite, then $G/N$ is an $\mathsf{MP}$-group.
\end{corollary}

In particular, if $G$ is a finitely generated nilpotent group with the
Magnus property then $G$ modulo its torsion subgroup $\tau(G)$ is a
finitely generated, torsion-free nilpotent group with the Magnus
property; in contrast, the finite group $\tau(G)$ does not in general
inherit the Magnus property; compare with
\cref{exa:fg-nilp-MP-any-class}.

The following sufficient criterion turns out to be useful for
rejecting the Magnus property and gives rise to the notion of basic
witness pairs which is to play a key role.

\begin{lemma} \label{lem:basic-witness-pair} Let $G$ be a group, and
  let $g \in G$ and
  $v \in [G,G] \smallsetminus \{ [g,w] \mid w \in G \}$ be such that
  $g^2 \not\equiv_{[G,G]} 1$ and
  $\langle g \rangle^G = \langle gv \rangle^G$.  Then $G$ is not an
  $\mathsf{MP}$-group.
\end{lemma}

\begin{proof}
  From $gv \equiv_{[G,G]} g \not \equiv_{[G,G]} g^{-1}$ it follows
  that $g$ and $gv$ are not inverse-conjugate in~$G$.
  They are not conjugate to one another either, as
  $gv \ne g [g,w] = g^w$ for all $w \in G$.  Thus
  $\langle g \rangle^G = \langle gv \rangle^G$ shows that $G$ does
  not have the Magnus property.
\end{proof}

For short we say that $(g,v)$ is a \emph{basic
  $\neg (\mathsf{MP})$-witness pair} for $G$, viz.\ a witness pair for
$G$ not having the Magnus property, if $g, v$ satisfy the conditions
in \cref{lem:basic-witness-pair}.  Part~(1) of the following lemma is
straightforward; compare with \cref{lem:retract}.  Part~(2) is
established by following the proof of \cref{prop:quot-crit} as per
contrapositive.

\begin{lemma} \label{lem:lift-witnesses}
  Let $G$ be a group, and let $N \nt G$.
  
  {\rm (1)} If $G = H \ltimes N$ splits over~$N$, then every basic
  $\neg (\mathsf{MP})$-witness pair for $H$ is also a basic
  $\neg (\mathsf{MP})$-witness pair for~$G$.

  {\rm (2)} Suppose that $N \subseteq [G,G]$ and that $g,v \in G$ are
  such that their images modulo $N$ form a basic
  $\neg (\mathsf{MP})$-witness pair $(\bar g, \bar v)$ for $G/N$.  If
  $\Omega_{gN} = \big\{ \langle gz \rangle^G \mid z \in N \big\}$
  satisfies the minimal condition, then $(\bar g, \bar v)$ lifts to a
  basic $\neg (\mathsf{MP})$-witness pair $(g_0,v_0)$ for $G$, with
  $g_0 \equiv_N g$ and $v_0 \equiv_N v$.
\end{lemma}

Another useful tool is the \textit{co-centraliser} of an element $g$
in a group $G$, that is
\[
  \cc_G(g) = \langle \, [g, w] \mid w \in G \rangle \le G;
\]
this group is closely related to the normal closure $\langle g
\rangle^G$ and thus of interest to us.

\begin{lemma}\label{lem:cocentraliser}
  Let $G$ be a group, and let $g \in G$.  Then
  \begin{enumerate}
  \item $\cc_G(g) \nt G$ and
    $\langle g \rangle^G = \langle g \rangle  \cc_G(g)$; in
    particular, $\langle g \rangle^G / \cc_G(g)$ is cyclic;
  \item $\cc_G(g) = [ \langle g \rangle^G,G]$, i.e., $\cc_G(g)$ is the
    smallest normal subgroup $N$ of $G$ such that
    $N \subseteq \langle g \rangle^G$ and $G$ acts trivially by
    conjugation on $\langle g \rangle^G/N$;
  \item $\cc_G(h) \leq \cc_G(g)$ for all $h \in \langle g \rangle^G$;
  \item if $h \in G$ with $\langle g \rangle^G = \langle h \rangle^G$
    then $\cc_G(g) = \cc_G(h)$.
  \end{enumerate}
\end{lemma}

\begin{proof}
  We set $X = \{ [g,w] \mid w \in G \}$ so that
  $\cc_G(g) = \langle X \rangle$ and
  $\langle g \rangle^G = \langle \{g\} \cup X \rangle$.  From the
  identity $[g,w]^v = [g,v]^{-1}[g,wv]$, for $w,v \in G$, we deduce
  that $\cc_G(g) \nt G$.  This establishes (1) and~(2).  Claims (3)
  and~(4) are immediate consequences of~(2).
\end{proof}

%%%
\section{Locally nilpotent groups}\label{sec:free_nilpotent_groups}

Clearly, if $G$ is an $\mathsf{MP}$-group then so is its
centre~$\Zent(G)$.  We are interested in sufficient conditions so that
the factor group $G/\Zent(G)$ inherits the Magnus property.  First we
observe a useful feature of co-centralisers in torsion-free, locally
nilpotent groups.

\begin{lemma}\label{lem:nilp->good}
  Let $G$ be a torsion-free, locally nilpotent group.  Let
  $g \in G\smallsetminus \Zent(G)$.  Then the co-centraliser of $g$
  satisfies
  \begin{equation}\label{eq:centre-in-co-centr}
    \langle g \rangle^G \cap \Zent(G) \leq \cc_G(g).
  \end{equation}
\end{lemma}

\begin{proof}
  Choose $v \in G$ such that $[g,v] \ne 1$.  For every
  $z \in \langle g \rangle^G \cap \Zent(G)$ there exist finitely many
  elements $w_1, \ldots, w_n \in G$ such that
  $z \in \langle g^{w_1}, \ldots, g^{w_n} \rangle$.  If the claim
  holds true for the nilpotent group
  $H = \langle g, v, w_1, \ldots, w_n \rangle \le G$, we conclude that
  $z \in \cc_H(g) \subseteq \cc_G(g)$.  Thus we may assume without
  loss that $G$ is nilpotent.
  
  Let $c$ denote the nilpotency class of~$G$, and let us fix the
  position where $g$ makes its appearance within the upper central
  series: $g \in \Zent_{i+1}(G)\smallsetminus \Zent_i(G)$ for suitable
  $i \in \{1, \dots, c-1\}$.  Since $G$ is torsion-free, so is
  $G/\Zent_i(G)$.  Hence we deduce from
  $\cc_G(g) \leq [\Zent_{i+1}(G),G] \leq \Zent_i(G)$ that
  $\langle g \rangle^G = \langle g \rangle \ltimes \cc_G(g)$ and
  consequently $\cc_G(g) = \langle g \rangle^G \cap \Zent_i(G)$.  This
  implies $\langle g \rangle^G \cap \Zent(G) \leq \cc_G(g)$.
\end{proof}

With this insight we are ready to deal with torsion-free, class-$2$
nilpotent groups.

\begin{proposition}\label{prop:class-2-nilp->mp}
  Let $G$ be a torsion-free, class-$2$ nilpotent group. Then $G$ has
  the Magnus property.
\end{proposition}

\begin{proof}
  Let $g, h \in G$ such that
  $\langle g \rangle^G = \langle h \rangle^G$.  If $g \in \Zent(G)$,
  then
  $\langle g \rangle = \langle g \rangle^G = \langle h \rangle^G =
  \langle h \rangle$ and, because $G$ is torsion-free, we conclude
  that $g \in \{h, h^{-1}\}$.

  Now suppose that $g \notin \Zent(G)$.  Since $G/\Zent(G)$ is
  torsion-free abelian, we deduce from
  $\langle g \rangle \Zent(G) = \langle h \rangle \Zent(G)$ that
  $g \equiv_{\Zent(G)} h^{\pm 1}$; replacing $g$ by its inverse if
  necessary, we may assume without loss that $g \equiv_{\Zent(G)} h$,
  hence $g^{-1} h \in \Zent(G) \cap \langle g \rangle^G \leq \cc_G(g)$
  by~\cref{lem:nilp->good}.  This implies
  $g^{-1} h = \prod_{i = 1}^k [g, v_i]^{e_i}$ for suitable
  $k \in \N_0$, $v_1, \ldots, v_k \in G$ and
  $e_1, \ldots, e_k \in \{1,-1\}$.  Since $G$ has nilpotency
  class~$2$, we obtain
  \[
    h = g \, \prod\nolimits_{i = 1}^k [g, v_i]^{e_i} = g \, \big[ g,
    \prod\nolimits_{i = 1}^{k} v_i^{\, e_i} \big] = g^v \qquad
    \text{for $v = \prod\nolimits_{i = 1}^{k} v_i^{\, e_i}$.} \qedhere
  \]
\end{proof}

\begin{example}
  It would be interesting to complement \cref{prop:class-2-nilp->mp}
  by characterising (or even classifying) finite, class-$2$ nilpotent
  groups with the Magnus property.  Each such group is necessarily a
  $\{2,3\}$-group and hence a direct product $G = P \times Q$ of its
  Sylow-$2$ and its Sylow-$3$ subgroup, each of which is again an
  $\mathsf{MP}$-group.

  However, it does not seem easy to give a succinct characterisation
  of finite, class-$2$ nilpotent $2$- or $3$-groups, in terms of
  canonical subgroups or quotients.  A halfway practical criterion
  for $3$-groups is the following: a finite,
  class-$2$ nilpotent $3$-group $G$ has the Magnus property if and
  only if (i) $\Zent(G)$ and $G/\Zent(G)$ are elementary abelian and
  (ii) for every ($\Zent(G)$-coset of an) element $g$ of order $9$
  there exists (a $\Zent(G)$-coset of) an element $h$ such that
  $[g,h] = g^3$.

  With this criterion it is, for instance, easy to see that the
  class-$2$ nilpotent group
  \[
    G = \langle t, a,b \mid t^3 = a^9 = b^9 = [a,b] = [a,t] b^3 =
    [b,t] a^3 = 1 \rangle,
  \]
  satisfies $\Zent(G) = [G,G] \cong C_3^{\, 2}$ and
  $G/\Zent(G) \cong C_3^{\, 3}$, but does not have the Magnus
  property; for instance, $a$ and $a^4$ have the same normal closure
  $\langle a, b^3 \rangle$, but are neither conjugate nor
  inverse-conjugate to one another.  Incidentally, examples of such
  kind illustrate that the condition of torsion-freeness is not
  redundant in \cref{lem:nilp->good}.
\end{example}

\begin{lemma}\label{lem:good->chains-of-length-1}
  Let $G$ be a group and let $g \in G \smallsetminus \Zent(G)$ be such
  that \cref{eq:centre-in-co-centr} holds.  Then every
  $z \in \langle g \rangle^G \cap \Zent(G)$ satisfies
  $\langle gz \rangle^G = \langle g \rangle^G$.
\end{lemma}

\begin{proof}
  Let $z \in \langle g \rangle^G \cap \Zent(G)$.  Clearly,
  $gz \in \langle g \rangle^G$ and it remains to show that
  $g \in \langle gz \rangle^G$.  Since
  $z \in \langle g \rangle^G \cap \Zent(G) \subseteq \cc_G(g)$, there
  exist $k \in \N$, $v_1, \dots, v_{k} \in G$ and
  $e_1, \ldots, e_k \in \{1,-1\}$ such that
  $z = \prod_{i = 1}^{k} [g, v_i]^{e_i}$.  Since $z$ is central, we
  deduce that
  $z = \prod_{i = 1}^{k} [gz, v_i]^{e_i} \in \cc_G(gz) \subseteq
  \langle gz \rangle^G$, and consequently
  $g \in \langle gz \rangle^G$.
\end{proof}

\begin{lemma}\label{lem:good->quotient}
  Let $G$ be an $\mathsf{MP}$-group, and suppose that
  \cref{eq:centre-in-co-centr} holds for each
  $g \in G \smallsetminus \Zent(G)$.  Then $G/\Zent(G)$ is an
  $\mathsf{MP}$-group.
\end{lemma}

\begin{proof}
  \cref{lem:good->chains-of-length-1} shows that, for each
  $g \in G\smallsetminus \Zent(G)$, any two distinct elements of
  $\Omega_{g\Zent(G)}$ are $\subseteq$-incomparable.  Thus
  \cref{prop:quot-crit} applies.
\end{proof}

From \cref{lem:nilp->good} and \cref{lem:good->quotient} we see that
within the class of torsion-free, locally nilpotent groups the Magnus
property passes from $G$ to $G/\Zent(G)$; this is a useful insight for
a future characterisation (or even classification) of finitely
generated, torsion-free nilpotent groups with the Magnus property.

\begin{corollary}\label{cor:reduction}
  Let $G$ be a torsion-free, locally nilpotent group.  If $G$ is an
  $\mathsf{MP}$-group then so is $G/\Zent(G)$.
\end{corollary}

For the next result we recall the notion of a basic
$\neg (\mathsf{MP})$-witness pair in the wake of
\cref{lem:basic-witness-pair}.

\begin{proposition}\label{prop:free-class>2-not-MP}
  Let $G$ be a free class-$c$ nilpotent group of rank at least~$2$,
  where $c \in \N_{\ge 3}$.  Then there exists a basic
  $\neg (\mathsf{MP})$-witness pair for~$G$.
\end{proposition}

\begin{proof}
  The free class-$c$ nilpotent group of rank $2$ is a retract of~$G$;
  by part~(1) of \cref{lem:lift-witnesses} we may suppose that
  $G = \langle x, y \rangle$ is freely generated by two elements.
  Furthermore, by part~(2) of \cref{lem:lift-witnesses} and induction
  on $c$ we may suppose that $c=3$.  In accordance with Witt's
  formula, the non-trivial sections of the lower central series of $G$
  are:
  $G / \gamma_2(G) = \langle \overline{x}, \overline{y} \rangle \cong
  C_\infty \times C_\infty$,
  $\gamma_2(G) / \gamma_3(G) = \langle \overline{[y,x]} \rangle \cong
  C_\infty$ and
  \[
    \Zent(G) = \gamma_3(G) = \langle [y,x,x], [y,x,y ] \rangle
    \cong C_\infty \times C_\infty;
  \]
  compare with~\cite[Chap.~3]{ClMaZy17}.

  Put $v = [y,x,y] \in \Zent(G)$.  Clearly, $x$ and $xv$ have the same
  normal closure in~$G$, namely
  $\langle x \rangle^G = \langle x \rangle \gamma_2(G) = \langle xv
  \rangle^G $.  Moreover, $x$ has infinite order modulo~$[G,G]$.  It
  remains to prove that $[x,w] \ne v$ for all~$w \in G$.

  Let $w \in G$.  For $[x,w] = v$ it is necessary that
  $[x,w] \in \Zent(G)$ and consequently
  $w \in \langle x \rangle \gamma_2(G)$.  But
  $w \equiv_{\Zent(G)} x^m [y,x]^n$ for $m,n \in \Z$ gives
  \[
    [x,w] = [x,[y,x]^n] = [y,x,x]^{-n} \ne [y,x,y] = v. \qedhere
  \]
\end{proof}

The following straightforward example illustrates that there are
finitely generated, nilpotent $\mathsf{MP}$-groups (with non-trivial
$3$-torsion) of any prescribed nilpotency class.

\begin{example} \label{exa:fg-nilp-MP-any-class}
  For $c \in \N$ the $2$-generated group
  \[
    G = \langle t,a \mid a^{3^c} = [a,t] a^{-3} = 1 \rangle \cong
    C_\infty \ltimes C_{3^c}
  \]
  is nilpotent of class $c$ and an
  $\mathsf{MP}$-group.

  Indeed, let $g,h \in G$ with
  $\langle g \rangle^G = \langle h \rangle^G$.  If $g=h=1$ then there
  is nothing to show.  Now suppose that $g$ and $h$ are non-trivial.
  There are unique parameters $l,m \in \Z$ with $0 \le m < 3^c$ such
  that $g = t^l a^m$.  Put
  $n = n(g) = 1+\min\{ v_3(l), v_3(m) \} \in \N$, where $v_3(k)$
  denotes the $3$-adic valuation of an integer~$k$.
  \cref{lem:cocentraliser} and a routine calculation show that
  \[
    M = \cc_G(g) = \cc_G(h) = \langle a^{3^n} \rangle = \{ [g,w] \mid w \in
    G \};
  \]
  further details are given at the end of
  Section~\ref{sec:tf-nilp-MP-gps}, where a related group $H$ is
  considered.  It suffices to show that $g \equiv_M h$ or
  $g \equiv_M h^{-1}$.

  If $l = 0$, then $n = 1 + v_3(m)$ and
  $\langle g \rangle^G = \langle g \rangle = \langle a^{3^{n-1}}
  \rangle$ gives
  $\langle g \rangle^G / M = \langle h \rangle^G / M \cong C_3$; it
  follows that $g \equiv_M h$ or $g \equiv_M h^{-1}$.  If $l \ne 0$,
  then $g$ and $h$ generate the same infinite cyclic subgroup
  modulo~$M$, and thus $g \equiv_M h$ or $g \equiv_M h^{-1}$.
\end{example}

It would be interesting to construct finitely generated, nilpotent
$\mathsf{MP}$-groups of prescribed nilpotency class that are
torsion-free.  In fact, it is already a challenge to construct an
explicitly a finitely generated, torsion-free, class-$3$ nilpotent
$\mathsf{MP}$-group, as we do below.  In principle,
\cref{cor:reduction} suggests that one proceeds by induction on the
nilpotency class.

\begin{lemma}\label{lem:class-3}
  Let $G$ be a group such that $\Zent(G)$ and $G/\Zent(G)$ are
  $\mathsf{MP}$-groups.  Suppose that for every
  $g \in G \smallsetminus \Zent(G)$, the set
  $\{ [g,w] \mid w \in G \}$ contains
  $\langle g \rangle^G \cap \Zent(G)$.  Then $G$ is an
  $\mathsf{MP}$-group.
\end{lemma}

\begin{proof}
  Suppose that $g, h \in G$ have the same normal closure in~$G$.  If
  $g \in \Zent(G)$ then $h \in \Zent(G)$; furthermore
  $\langle g \rangle = \langle h \rangle$ implies
  $g \in \{h, h^{-1}\}$, because $\Zent(G)$ is an $\mathsf{MP}$-group.

  Now suppose that $g \not\in \Zent(G)$.  Since $G/\Zent(G)$ is an
  $\mathsf{MP}$-group, there exist $v \in G$ and $e \in \{-1, 1\}$
  such that $g^v \equiv_{\Zent(G)} h^e$, and hence $g^v z = h^e$ for
  suitable~$z \in \langle g \rangle^G \cap \Zent(G)$.  Choose
  $w \in G$ such that $g^w = gz$; then
  $g^{wv} = (gz)^v = g^v z = h^e$.
\end{proof}

\begin{corollary} \label{cor:class-3} Let $G$ be a torsion-free group
  such that $G/\Zent(G)$ is an $\mathsf{MP}$-group.  If
  $\cc_G(g) \cap \Zent(G) \subseteq \{ [g,w] \mid w \in G \}$ for
  every $g \in G \smallsetminus \Zent(G)$, then $G$ is an
  $\mathsf{MP}$-group.
\end{corollary}

\begin{proof}
  Since $G$ is torsion-free, $\Zent(G)$ is an $\mathsf{MP}$-group.
  Let $g \in G \smallsetminus \Zent(G)$.  Since $G/\Zent(G)$ is
  torsion-free, we deduce that
  $\langle g \rangle^G \cap \Zent(G) = \cc_G(G) \cap \Zent(G)$.  Thus
  the claim follows from \cref{lem:class-3}.
\end{proof}

\begin{example} \label{exa:4-gen-class-3-MP-example}
  We use \cref{cor:class-3} to construct a $4$-generated,
  torsion-free, class-$3$ nilpotent $\mathsf{MP}$-group of Hirsch
  length~$9$.

  Let $F = \langle \dot x, \dot y, \dot z, \dot w \rangle$ be a free
  class-$3$ nilpotent group on four generators and consider
  \[
    G = \langle x,y,z,w \rangle = F \,/\, \langle \{ [\dot z, \dot y],
    [\dot w, \dot z] \} \cup R \rangle^F,
  \]
  where $x,y,z,w$ denote the images of
  $\dot x, \dot y, \dot z, \dot w$ and
  \[
    R = \left\{
      \begin{array}{llll}
      [\dot y,\dot x,\dot x], & [\dot y, \dot x, \dot y], & [\dot y, \dot x, \dot z] [\dot z,\dot x,\dot y]^{-1}, & [\dot y, \dot x, \dot w], \\
      {[\dot z, \dot x, \dot x],} & & [\dot z, \dot x, \dot z], &  [\dot z, \dot x, \dot w], \\
      {[\dot w, \dot x, \dot x]} [\dot z,\dot x,\dot y]^{-1}, & [\dot w,\dot x,\dot y], & [\dot w, \dot x, \dot z], & [\dot w, \dot x,\dot w], \\
       & [\dot w, \dot y, \dot y], & [\dot w,\dot y,\dot z] & [\dot w, \dot y, \dot w] [\dot z,\dot x,\dot y]^{-1}
      \end{array}
      \right\}.
  \]

  Clearly, $G$ is a $4$-generated nilpotent group of class at
  most~$3$.  The precise structure of $G$ can be determined as
  follows.  The collection process, subject to the initial ordering
  $\dot x < \dot y < \dot z < \dot w$, yields a Hall basis for $F$
  consisting of $4+6+20=30$ basic commutators; for instance,
  see~\cite[Section~3.1.3]{ClMaZy17}.  The relators in $R$ simply tell
  us to cancel or identify certain basic commutators of degree~$3$.
  The additional relators $[\dot z,\dot y]$ and $[\dot w,\dot z]$ are
  basic commutators of degree~$2$ to be cancelled; they also tell us
  to cancel the basic commutators
  $[\dot z, \dot y, \dot y], [\dot z, \dot y, \dot z], [\dot z, \dot
  y, \dot w], [\dot w, \dot z, \dot z], [\dot w, \dot z, \dot w]$ of
  degree~$3$.  The relations $[z,y,x]=[w,z,x]=[w,z,y]=1$, which also
  come with the relators $[\dot z,\dot y]$ and $[\dot w,\dot z]$, are
  already consequences of the relators in~$R$ and the Witt identity.
  For instance,
  $[z,y,x] = [z,y,x] [y,x,z] [z,x,y]^{-1} = [z,y,x] [y,x,z] [x,z,y] =
  1$.  In this way we see that $G$ is torsion-free of nilpotency
  class~$3$ and admits a poly-$C_\infty$ basis
  \[
    x,\; y,\; z,\; w, \quad [y, x],\; [z, x],\; [w, x],\; [w, y],\quad
    [y,x,z] = [z, x, y] = [w,x,x] = [w,y,w] 
  \]
  such that
  \[
    G / \gamma_2(G)= \langle \overline{x}, \overline{y}, \overline{z},
    \overline{w} \rangle \cong C_\infty^{\, 4}, \quad \gamma_2(G) /
    \gamma_3(G) = \langle \overline{[y, x]}, \overline{[z, x]},
    \overline{[w, x]}, \overline{[w, y]} \rangle \cong C_\infty^{\, 4}
  \]
  and
  \[
    Z = \Zent(G) = \gamma_3(G) = \langle [z,x,y] \rangle \cong
    C_\infty.
  \]
  In particular, $G$ has Hirsch length~$9$.  The group commutator
  induces a bi-additive map
  $\beta \colon \gamma_2(G)/Z \times G/\gamma_2(G) \to Z$ whose values
  on pairs of basis elements are given by the following commutator
  table:
  \[
    \begin{array}{c|llll}
      [\cdot,\cdot] & \; x & \; y & \; z & \; w \\\hline
      [y,x] & [y,x,x]=1 & [y, x, y]=1 & [y, x, z] = [z,x,y] & [y, x, w]=1 \\
      {[z,x]} & [z, x, x]=1 & [z,x,y] & [z, x, z]=1 &  [z, x, w]=1 \\
      {[w,x]} & [w, x, x] = [z,x,y] & [w,x,y]=1 & [w, x, z]=1 & [w, x,w]=1 \\
      {[w,y]} & \underline{[w, y, x]=1} & [w, y, y]=1 & [w,y,z]=1 &
                                                                    [w, y, w] = [z,x,y] \\ 
    \end{array}.
  \]
  The underlined entry is the only one that perhaps still requires a
  short explanation:
  \[
    [w,y,x] = [w,y,x][y,x,w][w,x,y]^{-1} = [w,y,x][y,x,w][x,w,y] = 1,
  \]
  using again relators from $R$ and the Witt identity.  The table
  shows that $\beta$ is a perfect pairing between $\gamma_2(G)/Z$ and
  $G/\gamma_2(G)$.
  
  It remains to verify that $G$ is an $\mathsf{MP}$-group.  By
  \cref{prop:class-2-nilp->mp}, the quotient $G/Z$ has the Magnus
  property.  By \cref{cor:class-3} it suffices to check that
  $\cc_G(g) \cap Z \subseteq \{ [g,v] \mid v \in G \}$ for all
  $g \in G \smallsetminus Z$.

  If $g \in \gamma_2(G)$ then
  $g \equiv_Z [y,x]^{m_1} [z,x]^{m_2} [w,x]^{m_3} [w,y]^{m_4}$ for
  suitable $m_1, \ldots, m_4 \in \Z$; since $\beta$ is a perfect
  pairing, this implies that
  $\cc_G(g) = \langle [z,x,y]^n \rangle = \{ [g,v] \mid v \in G \}$
  for $n = \gcd (m_1,m_2,m_3,m_4)$.  Now suppose that
  $g \not\in \gamma_2(G)$.  Since $G$ is nilpotent of class~$3$, the
  commutator identities
  \[
    [g,v]^{-1} = [g,v^{-1}] \underbrace{[g,v,v^{-1}]}_{\in Z}
    \quad \text{and} \quad [g,v] [g,w] = [g,wv]
    \underbrace{[g,w,v]^{-1}}_{\in Z}, \quad \text{for
      $v,w \in G$,}
  \]
  hold.  From these we deduce that every
  $h \in \cc_G(g) \cap Z$ is of the form $ h = h_1 h_2$ with
  \[
    h_1 \in \{ [g,v] \mid v \in G \} \cap Z \quad \text{and} \quad h_2
    \in \langle [g,w_1,w_2] \mid w_1,w_2 \in G \rangle \le Z.
  \]
  Observe that
  $g \equiv_{\gamma_2(G)} x^{m_1} y^{m_2} z^{m_3} w^{m_4}$, for
  suitable $m_1, \ldots, m_4 \in \Z$, and put
  $n = \gcd(m_1,m_2,m_3,m_4)$.  Since $\beta$ is a perfect pairing, we
  deduce as previously that
  \begin{multline*}
    \langle [g,w_1,w_2] \mid w_1,w_2 \in G \rangle = \{ [g,w_1,w_2]
    \mid w_1,w_2 \in G \} = \langle [z,x,y]^n \rangle \\ = \{ [g,w] \mid
    w \in \gamma_2(G) \}.
  \end{multline*}
  Writing $h_1 = [g,v]$ and $h_2 = [g,w]$, we obtain
  $h = h_1 h_2 = [g,v] [g,w]^v = [g,wv]$.
\end{example}

%%%%%

\section{Polynilpotent groups}\label{sec:polynilpotent_groups}

In this section we prove
\cref{thm:free-polynilpotent,thm:free-centre-by-polynilpotent}.  To
achieve this, we show first that the restricted wreath product
$C_\infty \wr C_\infty$ does not have the Magnus property.

\begin{proposition}\label{prop:wreath-product}
  There is a basic $\neg (\mathsf{MP})$-witness pair for the group
  $C_\infty \wr C_\infty$.
\end{proposition}

\begin{proof}
  We realise the wreath product as the group
  $G = \langle t \rangle \ltimes A \cong C_\infty \wr C_\infty$, where
  $A = \Z[T^{\pm 1}] = \Z[T,T^{-1}]$ is written additively and the
  action of $t$ on $A$ by conjugation is given by multiplication
  by~$T$.  The commutator of elements $x = t^m a$ and $y = t^n b$,
  with $m,n \in \Z$ and $a,b \in A$, is
  \[
    [x,y] = [t^m a, t^n b] = -a -b \cdot T^m + a \cdot T^n + b = a
    \cdot (T^n-1) - b \cdot (T^m - 1).
  \]
  In particular, $[G,G] = (T-1) \Z[T^{\pm 1}] \le A$ is the equal to
  the ideal of the ring $\Z[T^{\pm 1}]$ generated by $(T-1)$.
  If $x = t^ma \in G \smallsetminus A$ then
  $\langle x \rangle^G = \langle x \rangle \ltimes \cc_G(x)$ and
  \begin{equation} \label{equ:desc-I_x}
    \cc_G(x) = I_x = a \cdot (T-1) \, \Z[T^{\pm 1}] + (T^m-1) \,
    \Z[T^{\pm 1}]
  \end{equation}
  is the ideal of the ring $\Z[T^{\pm 1}]$ generated by
  $a \cdot (T-1)$ and $T^m-1$.
  
  Choose a prime $p \ge 5$ and consider the ring of integers
  $\mathcal{O} = \Z[\zeta] \cong \Z[T] / \Phi_p \Z[T]$ of the
  $p$\textsuperscript{th} cyclotomic field, where $\zeta$ denotes a
  primitive $p$\textsuperscript{th} root of unity and $\Phi_p$ the
  $p$\textsuperscript{th} cyclotomic polynomial.  It is well known
  that $\mathcal{O} / (\zeta-1) \mathcal{O} \cong \mathbb{F}_p$ is a
  field with $p$ elements.  By the Dirichlet unit theorem, the
  torsion-free rank of the unit group $\mathcal{O}^\times$ is
  $(p-1)/2 - 1 \ge 1$.  Consider the $(p-1)$\textsuperscript{st} power
  $\nu$ of an element of infinite order, for instance, the power of a
  cyclotomic unit such as $\nu = (\zeta+1)^{p-1}$, and write
  $\nu = f(\zeta)$ for a suitable polynomial
  $f \in 1 + (T-1) \, \Z[T]$.  Since $\nu$ has infinite order
  in~$\mathcal{O}^\times$, we deduce that
  \begin{equation} \label{equ:unit-inf-order} f \not\equiv_{T^p-1}
    T^n \qquad \text{for all $n \in \Z$.}
  \end{equation}
  Furthermore, we observe that $\nu^{-1}$ can be written in a similar
  form: $\nu^{-1} = {\bar f}(\zeta)$ for suitable
  ${\bar f} \in 1+ (T-1) \, \Z[T]$.  The embedding of rings
  \[
    \Z[T^{\pm 1}] / (T^p-1) \;\hookrightarrow\; \Z[T] / (T-1) \times
    \Z[T] / \Phi_p \Z[T] \; \cong \;  \Z \times \mathcal{O}
  \]
  shows that
  \begin{equation} \label{equ:expl-inv-of} f \cdot {\bar f}
    \equiv_{T^p-1} 1.
  \end{equation}

  Now consider $g = t^p e \in G$, where $e \in A$ denotes the constant
  polynomial~$1$, and
  $v = f - 1 \in (T-1) \Z(T) = [G,G] \subseteq A$.  Clearly, $g$
  has infinite order modulo $A \supseteq [G,G]$.  Commutators $[g,w]$,
  with $w = t^n b \in G$ for $n \in \Z$ and $b \in \Z[T^{\pm 1}]$, are
  of the form
  \[
    (T^n-1) - b \cdot (T^p-1) \equiv_{T^p-1} T^n-1;
  \]
  thus \cref{equ:unit-inf-order} shows that
  $v \not\in \{ [g,w] \mid w \in G \}$.  It remains to prove that
  $\langle g \rangle^G = \langle gv \rangle^G$, and for this it is
  enough to show that $v = f-1 \in (T-1) \, \Z[T^{\pm 1}]$ lies in
  $I_g \cap I_{gv}$.  Indeed, the general description provided in
  \cref{equ:desc-I_x} yields directly
  \[
    I_g = (T-1) \, \Z[T^{\pm 1}] + (T^p-1) \, \Z[T^{\pm 1}] = (T-1)
    \, \Z[T^{\pm 1}],
  \]
  and \cref{equ:expl-inv-of} implies that
  \[
    I_{gv} = f \cdot (T-1) \, \Z[T^{\pm 1}] + (T^p-1) \, \Z[T^{\pm 1}] = (T-1)
    \, \Z[T^{\pm 1}]. \qedhere
  \]  
\end{proof}

\begin{example}
  It easy to produce explicit examples of elements $g$ and $v$ in the
  proof of \cref{prop:wreath-product}, for which the claims could be
  checked directly by computation.  For $p=5$ one may take
  $\nu = (1+\zeta)^4 = 1 + (\zeta-1) f_0(\zeta)$ with
  $f_0 = 3+2T-T^2-2T^3$ and
  $\nu^{-1} = 1 + (\zeta-1) {\bar f}_0 (\zeta)$ with
  ${\bar f}_0 = -T +T^2-2T^3$.  A routine calculation yields
  \begin{multline*}
    \big( 1 + (T-1) f_0(T) \big) \big( 1 + (T-1) {\bar f}_0(T) \big) \\
    = 1 + (T-1) (3+T-4T^3) + (T-1)^2 ( -3T + T^2 -3T^3 - 3T^4 + 4T^6
      ) \\
    \equiv_{T^5-1}  1 + (T-1) (3+T-4T^3) + (T-1)^2 ( 3 + 4T + 4T^2 )
    = 1.
  \end{multline*}
\end{example}

Next we would like to use \cref{prop:wreath-product} to prove that, for
instance, the free metabelian group~$G$ of rank~$2$ does not have the
Magnus property.  If $C_\infty \wr C_\infty$ was a retract of~$G$,
then \cref{lem:retract} would immediately yield the desired
conclusion.  But, in fact, $C_\infty \wr C_\infty$ is not a retract
of~$G$.  Assume, for a contradiction, that $G = H \ltimes K$ with
$H \cong C_\infty \wr C_\infty$.  Then $[G,G] = [H,H] \ltimes [K,G]$
implies
$C_\infty \times C_\infty \cong G/[G,G] \cong H/[H,H] \times K/[K,G]
\cong C_\infty \times C_\infty \times K/[K,G]$, hence $[K,G] = K$.
But $G$ is residually a finite nilpotent group; compare
with~\cite{Gr57}.  Thus there is a finite-index normal subgroup
$N \trianglelefteq G$ such that $G/N$ is nilpotent and
$K \not \subseteq N$.  This gives $1 \ne KN/N \trianglelefteq G/N$
with $[KN/N,M/N] = KN/N$, a contradiction.

In the proof of the next proposition, which deals more generally with
free abelian-by-(class-$c$ nilpotent) groups, we by-pass the obstacle
that $C_\infty \wr C_\infty$ is not a retract.

\begin{proposition}\label{prop:ab-by-nilp}
  Let $c \in \N$, and let $G$ be an $\mathcal{N}_{(c,1)}$-free group
  of rank~$2$ viz.\ a free abelian-by-(class-$c$ nilpotent) group that
  is freely generated by two elements.  Then there exists a basic
  $\neg (\mathsf{MP})$-witness pair for~$G$.
\end{proposition}

\begin{proof}
  Let $H = \langle \bar x, \bar y \rangle = G/\gamma_{c+1}(G)$, the
  free class-$c$ nilpotent group of rank~$2$, and let $R = \Z H$
  denote the associated integral group ring.  Let
  $V = \mathbf{e} R \oplus \mathbf{f} R$ be the free right $R$-module
  on two free generators.  The Magnus embedding
  \[
    G \hookrightarrow H \ltimes V, \qquad
    w \mapsto \begin{pmatrix} \bar w & 0 \\ \mathbf{v}_w &
        1
      \end{pmatrix}
  \]
  allows us to realise the relatively free group $G$ as a group of
  $2 \times 2$ matrices with entries $\bar w \in H$ and
  $\mathbf{v}_w \in V$; compare with~\cite[\S~2.1]{Wi10}.  In this
  embedding, the matrices
  \[
    x = \begin{pmatrix} \bar x & 0 \\ \mathbf{e} & 1 \end{pmatrix}
    \qquad \text{and} \qquad y = \begin{pmatrix} \bar y & 0 \\
      \mathbf{f} & 1 \end{pmatrix}.
  \]
  constitute free generators of the relatively free group~$G$.

  Let us consider the iterated commutator
  \[
    z = [y,\, \underbrace{x, \ldots, x}_{c} \,] \in \gamma_{c+1}(G).
  \]
  We observe that~$z \ne 1$; for instance, one can deduce
  $z \not\in \gamma_{c+2}(G)$ from the fact that $z$ is a basic
  commutator in the Hall collection process, or carry out an explicit
  calculation as below.  It follows that the subgroup
  $\langle x, z \rangle = \langle x \rangle \ltimes A \le G$, with
  $\langle x \rangle \cong C_\infty$ and
  $A = \langle z^{x^m} \mid m \in \Z \rangle \le \gamma_{c+1}(G)$ free
  abelian, is isomorphic to the wreath product
  $C_\infty \wr C_\infty$: employing a transversal for the subgroup
  $\langle \bar x \rangle \cong C_\infty$ of~$H$, we can regard $V$ as
  a free $\Z \langle \bar x \rangle$-module (of infinite rank) and,
  accordingly, the non-trivial element $z$ of this module is moved
  about freely by the action of $x$ which is simply multiplication by
  the scalar $\bar x \in \Z \langle \bar x \rangle$.

  Following the proof of \cref{prop:wreath-product}, we consider
  elements
  \[
    g = x^p z \qquad \text{and} \qquad v = z^{f(x)-1},
  \]
  where $f \in 1+(T-1) \Z[T]$ is carefully chosen such that $g$ and
  $h = gv$ are not conjugate in
  $\langle x, z \rangle \cong C_\infty \wr C_\infty$, but generate the
  same normal closure in this subgroup and hence in~$G$.  Clearly,
  $g \equiv_{[G,G]} x^p \not\equiv_{[G,G]} 1$ has infinite order in
  $G/[G,G]$ and
  $v \in [\langle x, z \rangle,\langle x, z \rangle] \subseteq [G,G]$.
  It suffices to prove that $v \ne [g,w]$, or equivalently
  $g^w \ne h$, for all $w \in G$.

  As in the proof of \cref{prop:wreath-product}, let
  $\mathcal{O} = \Z[\zeta] \cong \Z[T]/\Phi_p \Z[T]$ denote the ring
  of integers of the $p$\textsuperscript{th} cyclotomic field, with
  $\zeta$ a primitive $p$\textsuperscript{th} root of unity; let
  $\nu = f(\zeta) \in \mathcal{O}^\times$ denote the one-unit of
  infinite order.  The kernel of the natural projection of rings
  $\pi \colon R \to \mathcal{O}$ specified by ${\bar x}^\pi = \zeta$
  and ${\bar y}^\pi = 1$ is generated, as a two-sided ideal, by the
  elements $\Phi_p(\bar x)$ and $\bar y -1$.  It induces a
  $\pi$-equivariant projection
  $\vartheta \colon V \to V/V \operatorname{ker}(\pi) \cong
  \dot{\mathbf{e}} \mathcal{O} \oplus \dot{\mathbf{f}} \mathcal{O}$,
  from the free $R$-module $V$ onto a free $\mathcal{O}$-module of
  rank~$2$ such that
  $(\mathbf{e} r + \mathbf{f} s)^\vartheta = \dot{\mathbf{e}} (r^\pi)
  + \dot{\mathbf{f}} (s^\pi)$ for all $r,s \in R$.
  
  It is straightforward to work out that
  \begin{align*}
    x^p=
    \begin{pmatrix} 
      {\bar x}^p & 0 \\
      \mathbf{v}_{x^p} & 1
    \end{pmatrix}, \quad
                       & \text{where ${\mathbf{v}_{x^p}}^\vartheta
                         = \dot{\mathbf{e}} \,\Phi_p(\zeta) = 0$,} \\
    z =
    \begin{pmatrix}
      1 & 0 \\
      \mathbf{v}_z & 1
    \end{pmatrix}, \quad
                     & \text{where ${\mathbf{v}_z}^\vartheta =
                       \dot{\mathbf{f}} (\zeta-1)^c$.} 
  \end{align*}
  From this we continue to see that
  \begin{align*}
    g  =
    \begin{pmatrix}
      {\bar x}^p & 0 \\
      \mathbf{v}_g & 1
    \end{pmatrix}, \quad
                       & \text{where ${\mathbf{v}_{g}}^\vartheta =
                         {\mathbf{v}_z}^\vartheta = \dot{\mathbf{f}}\, (\zeta-1)^c$,} \\  
    h =
    \begin{pmatrix}
      {\bar x}^p & 0 \\
      \mathbf{v}_h & 1
    \end{pmatrix}, \quad
                   & \text{where ${\mathbf{v}_h}^\vartheta =
                     ({\mathbf{v}_z}^\vartheta) \, f(\zeta) =
                     \dot{\mathbf{f}} \, (\zeta-1)^c \nu$.}
  \end{align*}
  Finally, we conjugate $g$ by an arbitrary element
  $w \in G$ to obtain
  \[
    \begin{pmatrix} {\bar g}^{\bar w} & 0 \\
      \mathbf{v}_{g^w} & 1 \end{pmatrix} = g^w =
    \begin{pmatrix} \bar w & 0 \\ \mathbf{v}_w &
      1 \end{pmatrix}^{-1} \begin{pmatrix} {\bar x}^p & 0 \\
      \mathbf{v}_g
      & 1 \end{pmatrix} \begin{pmatrix} \bar w & 0 \\
      \mathbf{v}_w & 1 \end{pmatrix}
    = \begin{pmatrix} {\bar w}^{-1} {\bar x}^p {\bar w} & 0 \\
      - \mathbf{v}_w {\bar w}^{-1} {\bar x}^p {\bar w} + \mathbf{v}_g
      {\bar w} + \mathbf{v}_w & 1 \end{pmatrix}
  \]
  and, for ${\bar w}^\pi = \zeta^m$ with suitable
  $m \in \{0,1,\ldots,p-1\}$, it follows that 
  \[
    ({\mathbf{v}_{g^w}})^\vartheta = ({\mathbf{v}_w}^\vartheta )
    \underbrace{(1-\zeta^p )}_{=0} + ( {\mathbf{v}_g}^\vartheta )
    \zeta^m = \dot{\mathbf{f}} \, (\zeta-1)^c \zeta^m.
  \]
  By construction, $\nu \in O^\times$ has infinite order so that
  $(\zeta-1)^c \nu \ne (\zeta-1)^c \zeta^m$.  Comparison with our
  computations for $h$ and~$g^w$ yields
  ${\mathbf{v}_h}^\vartheta \ne {\mathbf{v}_{g^w}}^\vartheta$, and
  hence $h \ne g^w$.
\end{proof}

We require another variant of
\cref{prop:free-class>2-not-MP,prop:ab-by-nilp}.  First, we record a
proposition which is perhaps folklore; we include a proof for
completeness.

\begin{proposition} \label{prop:small-centraliser} Let
  $G = \langle x,y \rangle$ be an $\mathcal{N}_\mathbf{c}$-free group
  of rank~$2$, where $l \in \N_{\ge 2}$ and $\mathbf{c} \in \N^l$.
  Let $n \in \N$.  Then the centraliser of $x^n$ in $G$ is
  $\operatorname{C}_G(x^n) = \langle x \rangle$.
\end{proposition}

\begin{proof}
  We write $\mathbf{c} = (c_1,\ldots,c_l)$,
  $\mathbf{c}' = (c_1,\ldots,c_{l-1})$ and put
  \[
    K = \gamma_{(c_1+1,\ldots,c_{l-1}+1)}(G) \nt G
  \]
  so that $G/K$ is an $\mathcal{N}_{\mathbf{c}'}$-free group of
  rank~$2$ and $K$ is a free class-$c_l$ nilpotent group.

  \smallskip
  
  \noindent \emph{Step 1.} First we argue by induction on $l$ that
  $\operatorname{C}_G(x^n) = \langle x \rangle
  \operatorname{C}_K(x^n)$.  Indeed, it suffices to fill in the base
  of the induction.  Suppose that $l=2$, and put $c = c_1$ so that
  $K = \gamma_{c+1}(G)$.  We observe that
  $H = \langle \bar x, \bar y \rangle = G/ \gamma_{c+2}(G)$ is a free
  class-$(c+1)$-nilpotent group of rank~$2$.  It suffices to show that
  $\operatorname{C}_H({\bar x}^n) = \langle \bar x \rangle
  \gamma_{c+1}(H)$.

  Clearly,
  $\operatorname{C}_H( {\bar x}^n) \subseteq \langle \bar x \rangle
  \gamma_2(H)$.  Hence it is enough to show that, for each
  $k \in \{2, \ldots,c\}$, the homomorphism of abelian groups
  \begin{align*}
    \gamma_k(H)/\gamma_{k+1}(H) &\to \gamma_{k+1}(H)/\gamma_{k+2}(H),\\
    w \gamma_{k+1}(H) &\mapsto [{\bar x}^n,w] \gamma_{k+2}(H) =
    [{\bar x},w]^n \gamma_{k+2}(H)
  \end{align*}
  is injective.  We may interpret the torsion-free sections
  $L_k = \gamma_k(H)/\gamma_{k+1}(H)$ as the first few homogeneous
  components of the free Lie ring $L = \bigoplus_{i=1}^\infty L_i$ on
  $\tilde x, \tilde y$, the images of $\bar x, \bar y$ in~$L_1$.
  Furthermore, we may think of $L$ as a Lie subring (generated by
  $\tilde x, \tilde y$) of the commutation Lie ring on a free
  associative ring~$A$, where $A$ is freely generated by non-commuting
  indeterminates $\tilde x, \tilde y$; compare
  with~\cite[Chap.~3]{ClMaZy17}.  The free ring~$A$ admits a natural
  $\N_0$-grading $A = \bigoplus_{i=0}^\infty A_i$, by means of the
  total $\{\tilde x, \tilde y\}$-degree function, which in turn
  induces the natural $\N$-grading $L = \bigoplus_{i=1}^\infty L_i$
  of~$L$.

  Fix $i \in \N_{\ge 2}$ and let $a \in L_i = L \cap A_i$ be a
  non-zero homogeneous Lie element of degree~$i$.  We are to show that
  the Lie commutator $[\tilde x, a]_\mathrm{Lie}$ is non-zero.  The
  monomials in $\tilde x, \tilde y$ of degree $i$ form a $\Z$-basis for
  the component $A_i$; we order them lexicographically
  \[
    {\tilde x}^i < {\tilde x}^{i-1} {\tilde y} < {\tilde x}^{i-2} {\tilde y}
    {\tilde x} < {\tilde x}^{i-2} {\tilde y}^2 < \ldots < {\tilde y}^{i-1}
    {\tilde x} < {\tilde y}^i,
  \]
  and proceed similarly with the monomials of degree~$i+1$.  Suppose
  that $a$, expressed as a $\Z$-linear combination of monomials, has
  leading term $m \, u(\tilde x, \tilde y)$ for
  $m \in \Z \smallsetminus \{0\}$ and $u(\tilde x, \tilde y)$ the smallest
  monomial occurring with non-zero coefficient.  Since $a \in L$ is a
  Lie element, we deduce that $u(\tilde x, \tilde y) \ne {\tilde x}^i$ and
  hence the Lie commutator
  $[\tilde x,a]_\mathrm{Lie} = {\tilde x} a - a {\tilde x} \in L_{i+1}
  \subseteq A_{i+1}$ is non-zero with leading term
  $m \, {\tilde x} u(\tilde x, \tilde y)$.
  
  \smallskip
  
  \noindent \emph{Step 2.}  It remains to prove that
  $\operatorname{C}_K(x^n) = 1$.  Put $c = c_l$.  Let
  $L = \bigoplus_{k=1}^c L_k$ denote the free class-$c$ nilpotent Lie
  ring associated to~$K$ and its lower central series; thus
  $L_k \cong \gamma_k(K)/\gamma_{k+1}(K)$ for $1 \le k \le c$ as a
  free $\Z$-module, and the Lie commutator of two homogeneous elements
  is induced by the group commutator, as in Step~$1$ above.  Extension
  of scalars yields the free class-$c$ nilpotent $\Q$-Lie algebra
  $\mathcal{L} = \bigoplus_{k=1}^c \mathcal{L}_k$, with
  $\mathcal{L}_k = \Q \otimes_\Z L_k$ for each~$k$.  Clearly,
  conjugation by $x$ induces an automorphism $\xi$ of~$\mathcal{L}$
  which respects the natural grading.  It suffices to prove that, for
  each~$k$, the only element of $\mathcal{L}_k$ fixed by $\xi^n$
  is~$0$.

  Put $H = G/K$ and $R = \Z H$.  The action of $G$ on $L$ factors
  through~$H$, and $\bar x \in H$ generates an infinite cyclic
  subgroup.  The Magnus embedding for the group $G/[K,K]$ shows that
  the $\Z \langle {\bar x} \rangle$-module $L_1$ embeds into a free
  $\Z \langle {\bar x} \rangle$-module (of infinite rank); compare
  with the proof of \cref{prop:ab-by-nilp}.  Thus, the
  $\Q \langle {\bar x} \rangle$-module $\mathcal{L}_1$ embeds into a
  free $\Q \langle {\bar x} \rangle$-module.  Observe that
  $\Q \langle {\bar x} \rangle$ is just the ring of Laurent
  polynomials over~$\Q$ and, in particular, a principal ideal domain.
  Therefore $\mathcal{L}_1$ is itself a free
  $\Q \langle {\bar x} \rangle$-module, with
  $\Q \langle {\bar x} \rangle$-basis
  $\mathbf{e}_1, \mathbf{e}_2, \ldots$, say.  Notice that
  $\mathcal{L}_1$ admits the $\Q$-basis
  \[
    \mathbf{f}_{i,m} = \mathbf{e}_i {\bar x}^m, \qquad
    \text{$i \in \N$ and $m \in \Z$;}
  \]
  these basis elements are at the same time free generators of  the
  free class-$c$ nilpotent $\Q$-Lie algebra~$\mathcal{L}$.

  Now fix $k \in \{1, \ldots, c\}$ and consider the action of $\xi$ on
  $\mathcal{L}_k$.  We observe that $\mathcal{L}_k$ is the $\Q$-span
  of the iterated Lie commutators
  \[
    \mathbf{F}_{\underline{i},\underline{m}} = [ \mathbf{f}_{i_1,m_1},
    \mathbf{f}_{i_2,m_2}, \ldots, \mathbf{f}_{i_k,m_k} ]_\mathrm{Lie},
  \]
  where $\underline{i} = (i_1, \ldots, i_k) \in \N^k$ and
  $\underline{m} = (m_1, \ldots, m_k) \in \Z^k$.  Furthermore, we
  understand the action of $\xi$ and hence of iterates $\xi^r$,
  $r \in \N$, on these Lie commutators:
  \[
    \mathbf{F}_{\underline{i},\underline{m}} \; \xi^r =
    \mathbf{F}_{\underline{i},\underline{m} +(r,r,\ldots,r)}.
  \]
  Let $\mathbf{v} \in \mathcal{L}_k \smallsetminus \{0\}$, written as a
  $\Q$-linear combination
  \[
    \mathbf{v} = \sum_{\underline{i} \in \N^k,\underline{m} \in \Z^k}
    v(\underline{i},\underline{m}) \,
    \mathbf{F}_{\underline{i},\underline{m}},
  \]
  where $v \colon \N^k \times \Z^k \to \Q$ is such that its `support'
  $S = \{ (\underline{i},\underline{m}) \in \N^k \times \Z^k \mid
  v(\underline{i},\underline{m}) \ne 0 \}$ is finite.  Also the `fine
  support' in the second coordinate
  \[
    S_\mathrm{fine} = \bigcup_{(\underline{i},\underline{m}) \in S} \{
    m_1, \ldots, m_k \} \subseteq \Z
  \]
  is finite.  Choose $r \in \N$ sufficiently large so that
  \[
    S_\mathrm{fine} \cap \big( S_\mathrm{fine} + rn \big) =
    \varnothing.
  \]
  Let $\mathcal{I} \nt \mathcal{L}$ be the Lie ideal generated by the
  following selection of free generators of~$\mathcal{L}$:
  $\mathbf{f}_{\underline{i},\underline{m}}$,
  $(\underline{i},\underline{m}) \in \N^k \times (\Z^k \smallsetminus
  S_\mathrm{fine})$.  Then
  $\mathbf{v} \not\equiv_\mathcal{I} 0 \equiv_\mathcal{I} \mathbf{v}
  \xi^{rn} = \mathbf{v} (\xi^n)^r$ implies
$\mathbf{v} \ne \mathbf{v} \xi^n$.
\end{proof}

\begin{proposition}\label{prop:nilp-cl-2-by-nilp}
  Let $c \in \N$, and let $G$ be an $\mathcal{N}_{(c,2)}$-free group
  of rank~$2$, viz.\ a free (class-$2$ nilpotent)-by-(class-$c$
  nilpotent) group that is freely generated by two elements.  Then
  there exists a basic $\neg (\mathsf{MP})$-witness pair for~$G$.
\end{proposition}

\begin{proof}
  The basic idea is to extend the proof of \cref{prop:ab-by-nilp}.
  For this purpose we put $K = \gamma_{c+1}(G)$, which is a free
  class-$2$ nilpotent group of countably infinite rank.  In
  particular, the abelianisation $K^\mathrm{ab} = K/[K,K]$ is free
  abelian, and $[K,K] \cong K^\mathrm{ab} \wedge K^\mathrm{ab}$ is the
  exterior square of $K^\mathrm{ab}$ and also free abelian.  Regarding
  $K^\mathrm{ab}$ as a free $\Z$-module, written additively, the
  exterior square is defined as
  \[
    K^\mathrm{ab} \wedge K^\mathrm{ab} = (K^\mathrm{ab} \otimes_\Z
    K^\mathrm{ab}) \;/\;
    \text{$\Z$-$\operatorname{span} \{ a \otimes b + b \otimes a \mid
      a,b \in K^\mathrm{ab} \}$}.
  \]
  Let $G = \langle x,y \rangle$, with free generators $x$ and $y$.
  The group commutators
  \[
    z_1 = [y,\, \underbrace{x, \ldots, x}_{c} \,] \in \gamma_{c+1}(G)
    \qquad \text{and} \qquad z_2 = [z_1,y] = [y,\, \underbrace{x, \ldots,
      x}_{c} \, , y] \in \gamma_{c+2}(G)
  \]
  lie in $K$ and their images $\breve z_1, \breve z_2$ modulo $[K,K]$
  yield $\Z$-linear independent generators of $K^\mathrm{ab}$; this
  follows, for instance, from the fact that $z_1$ and $z_2$ are basic
  commutators and form part of a Hall basis for
  $\gamma_{c+1}(G)/\gamma_{c+2}(G)$ and
  $\gamma_{c+2}(G)/\gamma_{c+3}(G)$, respectively, where
  $\gamma_{c+3}(G) \supseteq \gamma_{2c+2}(G) \supseteq [K,K]$.  Hence
  the group commutator
  \[
    z = [z_1,z_2] \in [K,K] \smallsetminus \{1\}
  \]
  is non-trivial; in the additive notation, it corresponds to
  $\breve z_1 \wedge \breve z_2 \in K^\mathrm{ab} \wedge K^\mathrm{ab}
  \smallsetminus \{0\}$.  The action of $G$ on $K^\mathrm{ab}$ and on
  $[K,K] \cong K^\mathrm{ab} \wedge K^\mathrm{ab}$ factors through
  $H = G/K$; concretely,
  $z^w = [z_1,z_2]^w = [z_1^{\, w}, z_2^{\, w}]$ translates to
  $(\breve z_1 \wedge \breve z_2).{\bar w} = ({\breve z_1}.{\bar w})
  \wedge ({\breve z_2}.{\bar w})$ for $w \in G$ with image
  $\bar w \in H$.
  
  Let $R = \Z H$ denote the integral group ring associated
  to~$H = \langle \bar x, \bar y \rangle$, and let
  $V = \mathbf{e} R \oplus \mathbf{f} R$ denote the free right
  $R$-module of rank~$2$.  We compose reduction modulo
  $[K,K]$ with the Magnus embedding for $G/[K,K]$ to obtain a
  homomorphism
  \[
    \eta \colon G \to G/[K,K] \hookrightarrow H \ltimes V, \qquad w
    \mapsto \begin{pmatrix} \bar w & 0 \\ \mathbf{v}_w & 1
      \end{pmatrix};
  \]
  compare with the proof of \cref{prop:ab-by-nilp}.  The generators
  $x,y$ of $G$ are mapped to
  \[
    x \eta = \begin{pmatrix} \bar x & 0 \\ \mathbf{e} &
      1 \end{pmatrix}
    \qquad \text{and} \qquad  y \eta = \begin{pmatrix} \bar y & 0 \\
      \mathbf{f} & 1 \end{pmatrix}.
  \]

  The exterior square $V \wedge V$ of the $\Z$-module $V$ is an
  $R$-module via the diagonal action.  In fact, $V \wedge V$ is a free
  $R$-module (of countably infinite rank): if
  $H_0 \subseteq H \smallsetminus \{1\}$ is a set of representatives
  for the equivalence classes
  $\{ \bar w, {\bar w}^{-1} \} \subseteq H \smallsetminus \{1\}$ of
  the relation ``equal or inverse to one another'', then the elements
  \[
    \mathbf{e} \wedge \mathbf{e} {\bar w}, \;\; \text{(for
      $\bar w \in H_0$)}, \quad \mathbf{f} \wedge \mathbf{f} {\bar w}
    \;\; \text{(for $\bar w \in H_0$)}, \quad \mathbf{e} \wedge
    \mathbf{f} {\bar w} \;\; \text{(for $\bar w \in H$)}
  \]
  constitute an $R$-basis for $V \wedge V$.  Since the $R$-module
  $K^\mathrm{ab}$ embeds into~$V$, the $R$-module
  $K^\mathrm{ab} \wedge K^\mathrm{ab}$ embeds into the free $R$-module
  $V \wedge V$.  A similar argument as in the proof of
  \cref{prop:ab-by-nilp} shows that the subgroup
  $\langle x, z \rangle = \langle x \rangle \ltimes \langle z^{x^m}
  \mid m \in \Z \rangle$ is isomorphic to the wreath product
  $C_\infty \wr C_\infty$.  As before, we consider elements
  \[
    g = x^p z \qquad \text{and} \qquad v = z^{f(x)-1},
  \]
  where $f \in 1+(T-1) \Z[T]$ is such that $g$ and $gv$ are not
  conjugate in $\langle x, z \rangle$, but generate the same normal
  closure in this subgroup and hence in~$G$.  Clearly,
  $g \equiv_{[G,G]} x^p \not\equiv_{[G,G]} 1$ has infinite order in
  $G/[G,G]$ and $v \in [K,K] \subseteq [G,G]$.  It suffices to prove
  that $[g,w] \ne v$ for all $w \in G$.  Reduction modulo $[K,K]$
  shows that for $[g,w] = v$ it would be necessary that
  $ [x^p,w] \equiv_{[K,K]} 1$; hence by \cref{prop:small-centraliser}
  we only need to prove that $[g,w] \ne v$ for $w \in [K,K]$.

  As before let $\mathcal{O} = \Z[\zeta]$ denote the ring of integers
  of the $p$\textsuperscript{th} cyclotomic field, with $\zeta$ a
  primitive $p$\textsuperscript{th} root of unity; by construction,
  $\nu = f(\zeta) \in \mathcal{O}^\times$ has infinite order.  We
  consider the ring of Laurent polynomials
  $\mathcal{O}[Y^{\pm 1}] = \mathcal{O}[Y,Y^{-1}]$.  The natural
  projection of rings $\pi \colon R \to \mathcal{O}[Y^{\pm 1}]$
  specified by ${\bar x}^\pi = \zeta$ and ${\bar y}^\pi = Y$ induces a
  $\pi$-equivariant projection
  $\vartheta \colon V \to \dot{\mathbf{e}} \mathcal{O}[Y^{\pm 1}]
  \oplus \dot{\mathbf{f}} \mathcal{O}[Y^{\pm 1}] = \dot V$, from the
  free $R$-module $V$ onto a free $\mathcal{O}[Y^{\pm 1}]$-module
  $\dot V$ on $\dot{\mathbf{e}}, \dot{\mathbf{f}}$.
  
  It is straightforward to work out that
  \begin{align*}
    z_1 \eta = 
    \begin{pmatrix}
      1 & 0 \\
      \mathbf{v}_{z_1} & 1
    \end{pmatrix}, \quad
                     & \text{where ${\mathbf{v}_{z_1}}^\vartheta =
                       \dot{\mathbf{e}} \, (1-Y) (\zeta-1)^{c-1} +
                       \dot{\mathbf{f}} \, (\zeta-1)^c$.} \\
    z_2 \eta =
    \begin{pmatrix}
      1 & 0 \\
      \mathbf{v}_{z_2} & 1
    \end{pmatrix}, \quad
                     & \text{where ${\mathbf{v}_{z_2}}^\vartheta =
                       \dot{\mathbf{e}} \, (1-Y)^2 (\zeta-1)^{c-1} +
                       \dot{\mathbf{f}} \, (1-Y) (\zeta-1)^c$.}
  \end{align*}
  Restriction of scalars turns $\dot V$ into a free
  $\mathcal{O}$-module, with an $\mathcal{O}$-basis consisting of
  \[
    \dot{\mathbf{e}} Y^m \; \text{($m \in \Z$)}, \quad
    \dot{\mathbf{f}} Y^n \; \text{($n \in \Z$)}.
  \]
  Thus the exterior square $\dot V \wedge_\mathcal{O} \dot V$ over
  $\mathcal{O}$ is a free $\mathcal{O}$-module, with
  $\mathcal{O}$-basis
  \[
    \dot{\mathbf{e}} Y^m \wedge \dot{\mathbf{e}} Y^n \; \text{(for
      $m<n$)}, \;\; \dot{\mathbf{f}} Y^m \wedge \dot{\mathbf{f}} Y^n \; \text{(for
      $m<n$)}, \;\; \dot{\mathbf{e}} Y^m \wedge \dot{\mathbf{f}} Y^n, \quad
    \text{where $m,n \in \Z$;}
  \]
  below we express elements of $\dot V \wedge_\mathcal{O} \dot V$ with
  respect to this $\mathcal{O}$-basis, keeping track of the
  coefficients of the basis element
  $\dot{\mathbf{e}} Y \wedge \dot{\mathbf{e}} Y^2$.
  
  Next we consider the composition $\psi = \hat \vartheta \varrho$ of
  the $\pi$-equivariant morphism
  \[
    \hat \vartheta \colon [K,K] \cong K^\mathrm{ab} \wedge
    K^\mathrm{ab} \hookrightarrow V \wedge V \to \dot V \wedge \dot V
  \]
  with the canonical homomorphism of $\Z$-modules
  $\varrho \colon \dot V \wedge \dot V \to \dot V \wedge_{\mathcal{O}}
  \dot V$.  A routine computation yields
  \begin{multline*}
    z^\psi = [z_1,z_2]^{\hat \vartheta \varrho} = \big(
    {\mathbf{v}_{z_1}}^\vartheta
    \wedge {\mathbf{v}_{z_2}}^\vartheta \big)^\varrho \\
    = \Big( \Big( \big( \dot{\mathbf{e}} \, (1-Y) +
    \dot{\mathbf{f}} \, (\zeta-1) \big)  \wedge \big(
    \dot{\mathbf{e}} \, (1-Y)^2  + \dot{\mathbf{f}} \,
    (Y-1) (\zeta-1) \big) \Big) (\zeta-1)^{2c-2} \Big)^\varrho\\
    = \ldots + (\dot{\mathbf{e}} Y \wedge \dot{\mathbf{e}} Y^2) \,
    (-(\zeta-1)^{2c-2}) + \ldots,
  \end{multline*}
  where on the far right-hand side we only display the
  $\dot{\mathbf{e}} Y \wedge \dot{\mathbf{e}} Y^2$-term.  From this we
  deduce that
  \[
    v^\psi = \big( z^{f(x) -1} \big)^\psi = z^\psi \, (f(\zeta)-1)
    = \ldots + (\dot{\mathbf{e}} Y \wedge \dot{\mathbf{e}} Y^2) \,
    (-(\zeta-1)^{2c-2} (\nu -1)) + \ldots,
  \]
  where we again only record the
  $\dot{\mathbf{e}} Y \wedge \dot{\mathbf{e}} Y^2$-term to see that
  $v^\psi$ is non-zero.

  Finally, we deduce that $[g,w] \ne v$ for all $w \in [K,K]$ from
  \[
    [g,w]^\psi = [x^p z, w]^\psi = [x^p,w]^\psi = (w^{1-x^p})^\psi =
    w^\psi \, (1- \zeta^p) = 0 \ne v^\psi. \qedhere
  \]
\end{proof}

\begin{proof}[Proof of \cref{thm:free-polynilpotent}]
  The group $G$ is an $\mathcal{N}_\mathbf{c}$-free group of rank~$d$,
  for parameters $d,l \in \N$ and $\mathbf{c} \in \N^l$.  If $d=1$ or
  $\mathbf{c} \in \{ (1), (2) \}$ then $G$ is nilpotent of class at
  most $2$; in this situation \cref{prop:class-2-nilp->mp} implies
  that $G$ has the Magnus property.

  Now suppose that $d \ge 2$ and that
  $\mathbf{c} \not\in \{ (1), (2) \}$.  We are to show that $G$ does
  not have the Magnus property, and by \cref{cor:f-g-enough} we may
  suppose that $d=2$.  If $l = 1$ then $G$ is an
  $\mathcal{N}_{(c)}$-free group with $c = c_1 \ge 3$, and
  \cref{prop:free-class>2-not-MP} shows that $G$ does not have the
  Magnus property.  Likewise, if $l = 2$ and $c_2 \in \{1,2\}$, then
  $G$ is an $\mathcal{N}_{(c,1)}$-free or $\mathcal{N}_{(c,2)}$-free
  group of rank $2$ for $c = c_l$, and $G$ does not have the Magnus
  property by \cref{prop:ab-by-nilp,prop:nilp-cl-2-by-nilp}.

  Thus, we may suppose that we are in none of these special
  circumstances.  We write $\mathbf{c} = (c_1,\ldots,c_l) \in \N^l$
  and distinguish two cases.

  \smallskip

  \noindent \emph{Case~$1$}: $c_l \ge 3$.  In this case $l \ge 2$ and
  we put $N = \gamma_{(c_1+1,\ldots,c_{l-1}+1)}(G)$.  We note that
  $G/N$ is an $\mathcal{N}_{\mathbf{c}'}$-free group of rank $2$, for
  $\mathbf{c}' = (c_1,\ldots,c_{l-1})$, while $N$ is a free class-$c$
  nilpotent group with $c = c_l \ge 3$ of countably infinite rank.

  \smallskip

  \noindent \emph{Case~$2$}: $c_l \in \{1,2\}$.  In this case
  $l \ge 3$ and we put $N = \gamma_{(c_1+1,\ldots,c_{l-2}+1)}(G)$.  We
  note that $G/N$ is an $\mathcal{N}_{\mathbf{c}'}$-free group of rank
  $2$, for $\mathbf{c}' = (c_1,\ldots,c_{l-2})$, while $N$ is an
  $\mathcal{N}_{(c,1)}$-free or $\mathcal{N}_{(c,2)}$-free group with
  $c = c_{l-1}$ of countably infinite rank.

  \smallskip

  In any case,
  \cref{prop:free-class>2-not-MP,prop:ab-by-nilp,prop:nilp-cl-2-by-nilp}
  and part~(1) of \cref{lem:lift-witnesses} provide $g,v \in N$ such
  that $(g,v)$ is a basic $\neg (\mathsf{MP})$-witness pair for~$N$.
  Clearly, $g$ and $h=gv$ have the same normal closure
  $\langle g \rangle^G = \langle h \rangle^G$ in~$G$, and it suffices
  to prove that $g$ and $h$ are neither conjugate nor
  inverse-conjugate to one another in~$G$.

  For a contradiction, assume that $g^w \in \{ h, h^{-1} \}$ for some
  $w \in G$.  We put $H = G/N$ and $R = \Z H$.  Observe that $G/[N,N]$
  is a free abelian-by-$\mathcal{N}_{\mathbf{c}'}$ group of rank~$2$.
  The Magnus embedding for this group yields an embedding of the
  $R$-module $N^\mathrm{ab} = N/[N,N]$ into a free $R$-module. 
  %; compare with~\cite[\S~2.1]{Wi10}.  
  Our assumption yields
  $\mathbf{v}_g {\bar w} \in \{ \mathbf{v}_g, - \mathbf{v}_g \}$,
  hence $\mathbf{v}_g ({\bar w} -1)=0$ or
  $\mathbf{v}_g ({\bar w}+1) = 0$, where $\mathbf{v}_g$ denotes the
  image of $g$ in $N^\mathrm{ab}$, regarded as a module element, and
  ${\bar w} \in R$ denotes the image of $w$ in $H \subseteq R$.  We
  observe that $g \not\in [N,N]$ implies that $\mathbf{v}_g \ne 0$.
  The group $H$ is right-orderable and thus the group ring $R$ has no
  zero-divisors; see \cite[Ch.~13, Thm.~1.11]{Pa77}, where the result
  is attributed to Bovdi~\cite{Bo60}.  This implies ${\bar w}-1=0$ or
  ${\bar w}+1=0$ in~$R$.  From ${\bar w} \in H$ we see that
  ${\bar w} \ne -1$.  Hence ${\bar w} = 1$ and $w \in N$, in
  contradiction to the initial choice of $g$ and $h=gv$ which
  precludes that they are conjugate in~$N$.
\end{proof}

Finally, we extend \cref{thm:free-polynilpotent} to yet another class
of relatively free groups.  For any $d,l \in \N$ and
$\mathbf{c} = (c_1,\ldots,c_l) \in \N^l$, the \emph{free
  centre-by-$\mathcal{N}_\mathbf{c}$ group} of rank $d$ can be
constructed as the quotient $F/[\gamma_{(c_1+1,\ldots,c_l+1)},F]$ of
an absolutely free group of rank~$d$.

\begin{proposition} \label{prop:centre-by-polynilp-not}
  Let $G$ be a free centre-by-$\mathcal{N}_\mathbf{c}$ group of
  rank~$2$, where $\mathbf{c} \in \N^l$ with $l \in \N_{\ge 2}$.  Then
  there exists a basic $\neg (\mathsf{MP})$-witness pair for the
  group~$G$.
\end{proposition}

\begin{proof}
  Write $G = \langle x,y \rangle$, with free generators $x,y$, and
  $\mathbf{c} = (c_1,\ldots,c_l)$.  Let
  $Z = \gamma_{(c_1+1,\ldots,c_l+1)}(G) \subseteq [G,G]$.  Consider
  $g = x$ and any $v \in Z \smallsetminus \{1\}$.  Clearly, $g$ has
  infinite order modulo $[G,G]$.  It remains to show that
  $\langle g \rangle^G = \langle gv \rangle^G$ and
  $v \not\in \{ [g,w] \mid w \in G \}$.

  Since $v \in [G,G] \subseteq \langle x \rangle^G$, we find
  $k \in \N$, $e_1, \ldots, e_k \in \{1,-1\}$ and
  $w_1, \ldots, w_k \in G$ such that
  $v = \prod_{i=1}^k (x^{e_i})^{w_i}$.  From
  $x^{\sum_{i=1}^k e_i} \equiv_{[G,G]} v \equiv_{[G,G]} 1$ we conclude
  that $\sum_{i=1}^k e_i =0$.  Since $v$ is central in~$G$, this gives
  \[
    v = v^{\sum_{i=1}^k e_i} \prod\nolimits_{i=1}^k (x^{e_i})^{w_i} =
    \prod\nolimits_{i=1}^k \big( (xv)^{e_i} \big)^{w_i} \in \langle g
    \rangle^G \cap \langle gv \rangle^G
  \]
  and hence $\langle g \rangle^G = \langle gv \rangle^G$.

  Next assume, for a contradiction, that $[g,w] = v$ for some
  $w \in G$.  Then $[x,w] \equiv_Z 1$, and
  \cref{prop:small-centraliser} implies $w = x^m z$, for suitable
  $m \in \Z$ and $z \in Z$.  This gives $[g,w] = [x,x^mz] = 1 \ne v$,
  a contradiction.
\end{proof}

\begin{proof}[Proof of \cref{thm:free-centre-by-polynilpotent}]
  The group $G$ is a free centre-by-$\mathcal{N}_\mathbf{c}$ group $G$
  of rank~$d$, for parameters $d,l \in \N$ and $\mathbf{c} \in \N^l$.
  If $d=1$ or $\mathbf{c} = (1)$ then $G$ is nilpotent of class at
  most $2$; in this situation \cref{prop:class-2-nilp->mp} implies
  that $G$ has the Magnus property.

  Now suppose that $d \ge 2$ and that $\mathbf{c} \ne (1)$.  We are to
  show that $G$ does not have the Magnus property, and by
  \cref{cor:f-g-enough} we may suppose that $d=2$.  If $l = 1$ then
  $G$ is an $\mathcal{N}_{(c+1)}$-free group with $c = c_1 \ge 2$, and
  \cref{prop:free-class>2-not-MP} shows that $G$ does not have the
  Magnus property.  If $l \ge 2$ then
  \cref{prop:centre-by-polynilp-not} shows that $G$ does not have the
  Magnus property.
\end{proof}

%%%
\section{Torsion-free nilpotent groups with the Magnus
  property} \label{sec:tf-nilp-MP-gps}

In this section we establish \cref{thm:nilpotent-existence}.  Let
$c \in \N$.  In order to construct a torsion-free, nilpotent
$\mathsf{MP}$-group of prescribed nilpotency class~$c$ we aim to build
an ultraproduct
$\mathcal{G} = \big( \prod\nolimits_{p \in P} G_p \big) /
\!\!\sim_\mathfrak{U}$ of suitable finite $p$-groups $G_p$, where
$P = \mathbb{P}_{>2}$ denotes the set of all odd primes, and to appeal
to {\L}o{\'s}'s theorem.

Being class-$c$ nilpotent, not having $q$-torsion for a given
prime~$q$, and possessing the Magnus property are first-order
properties; thus it would suffice to construct a family of finite
nilpotent $\mathsf{MP}$-groups $G_p$, $p \in P$, such that each group
has nilpotency class $c$ and such that for any prime $q$ there are
only finitely many $p \in P$ with $q \mid \lvert G_p \rvert$.
However, finite nilpotent $\mathsf{MP}$-groups are necessarily
$\{2,3\}$-groups.  More generally, every group with the Magnus
property is \textit{inverse semi-rational}, that is, every pair of
elements generating the same subgroup (not necessarily normal) is
already a pair of conjugate or inverse-conjugate elements. Finite
groups with this property can be characterised using character
theory. Chillag and Dolfi \cite{ChDo10} establish that all finite
soluble inverse semi-rational groups are $\{2, 3, 5, 7, 13\}$-groups.
Consequently, no family as described above exists. However, we can
salvage our strategy by considering a variant of the Magnus property.

For convenience we use $[k,l]_\Z = \{ m \in \Z \mid k \le m \le l \}$,
for $k,l \in \Z$, as a short notation for intervals in~$\Z$.
Suppose that $G_p$, $p \in P$, is a family of groups with the
following properties:
\begin{enumerate}[\rm (i)]
\item for each~$p$, the group $G_p$ is a metabelian finite $p$-group
  of nilpotency class $c$;
\item there exists a non-decreasing function $f \colon \N \to \N$ such
  that, for each $p \in P$, the group $G_p$ satisfies the following
  uniform, but `weak' Magnus property:
  \begin{equation} \label{equ:magnus-f} \tag{$\mathsf{wM}_f$}
    \begin{split} 
      & \forall g,h \in G_p \quad \forall N \in \N \quad
      \forall k,l \in [0,N]_\Z \\
      & \qquad \forall e_1, \ldots, e_k \in \{1,-1\} \quad \forall
      v_1, \ldots, v_k \in G_p \\
      & \qquad \forall d_1, \ldots, d_l \in \{1,-1\} \quad \forall
      w_1, \ldots, w_l \in G_p : \\
      & \Big( h = \prod\nolimits_{i=1}^k (g^{e_i})^{v_i} \; \land \; g
      =
      \prod\nolimits_{j=1}^l (h^{d_j})^{w_j} \Big) \\
      & \; \implies \; \Big( \exists r,s \in [-f(N),f(N)]_\Z \quad
      \exists v,w \in G_p :\quad g = (h^r)^v \; \land \; h = (g^s)^w
      \Big).
    \end{split}
  \end{equation}
\end{enumerate}
We observe that the quantifier over the integer $N$ can be eliminated
by passing to a countable collection of sentences in the first-order
language of groups; the quantifiers over $k,l$ are purely for
convenience and can be eliminated directly.

Let $\mathfrak{U}$ be a non-principal ultrafilter on the index set
$P$. Then, by {\L}o{\'s}'s theorem, the ultraproduct
\[
  \mathcal{G} = \left( \prod\nolimits_{p \in P} G_p  \right)
  / \sim_\mathfrak{U} 
\]
is a metabelian, torsion-free, class-$c$ nilpotent group satisfying
the uniform `weak' Magnus property~\eqref{equ:magnus-f}; compare
with \cite[Thm.~4.1.9]{ChKe90}.  But, since $\mathcal{G}$ is torsion-free
nilpotent, the latter implies that $\mathcal{G}$ has the Magnus
property.  Indeed, suppose that $g,h \in \mathcal{G}$ are such that
$\langle g \rangle^G = \langle h \rangle^G$.  If $g = 1$ then $h = 1$,
and $g = h$ so that $g$ and $h$ are certainly conjugate.  Now suppose
that $g \ne 1$.  Then \eqref{equ:magnus-f} yields $r,s \in \Z$ and
$v,w \in \mathcal{G}$ such that $g = (h^r)^v$ and $h = (g^s)^w$, thus
\[
  g = (h^r)^v = \Big( \big( (g^s)^w \big)^r
  \Big)^v = (g^{rs})^{vw}.
\]
Consider the upper central series
$1 = \mathrm{Z}_0(\mathcal{G}) \le \mathrm{Z}_1(\mathcal{G}) \le
\ldots \le \mathrm{Z}_c(\mathcal{G}) = \mathcal{G}$ of the nilpotent
group~$\mathcal{G}$.  Since $g \ne 1$, we find $i \in [1,c]_\Z$ such
that
$g \in \mathrm{Z}_i(\mathcal{G}) \smallsetminus
\mathrm{Z}_{i-1}(\mathcal{G})$.  Since
$\mathrm{Z}_i(\mathcal{G})/\mathrm{Z}_{i-1}(\mathcal{G})$ is
torsion-free, $g$ generates an infinite cyclic group modulo
$\mathrm{Z}_{i-1}(\mathcal{G})$ and the congruence
$g \equiv_{\mathrm{Z}_{i-1}(\mathcal{G})} g^{rs}$ implies that
$rs = 1$. Thus $r \in \{1,-1\}$, and $g = (h^r)^v$ is conjugate to $h$
or to $h^{-1}$ in $\mathcal{G}$.

Finally, since all relevant properties of $\mathcal{G}$, including the
ordinary Magnus property, are expressible in terms of first-order
sentences, the L\"owenheim--Skolem theorem \cite[Cor.~2.1.4]{ChKe90}
shows that there exists a countable $\mathsf{MP}$-group
$\dot{\mathcal{G}}$ which is metabelian, torsion-free and nilpotent of
class precisely~$c$.  This establishes \cref{thm:nilpotent-existence}.

\medskip

It remains to construct the family of groups $G_p$, $p \in P$, with
the properties (i) and (ii) described above.  We give one rather
concrete construction.  Fix an odd prime~$p$, and consider
\begin{equation} \label{equ:concrete-Gp}
  G = G_p = \langle t, a \mid [a,t] = a^p, \; t^{p^{c-1}} = a^{p^c} =
  1 \rangle.
\end{equation}
Clearly, $G = \langle t \rangle \ltimes \langle a \rangle$ is
metacyclic, with $\langle t \rangle \cong C_{p^{c-1}}$ and
$\langle a \rangle \cong C_{p^c}$.  It is easy to work out the lower
central series:
\[
  \gamma_1(G) = G \qquad \text{and} \qquad \gamma_{i+1}(G) = \langle
  a^{p^i} \rangle \qquad \text{for $i \in \N$;}
\]
in particular, $G$ has nilpotency class~$c$.  In order to check the
`weak' Magnus property, we make use of the following lemma, which is
heuristically a torsion analogue of \cref{prop:class-2-nilp->mp}.

\begin{lemma} \label{lem:commutator-set-closed} Let $G$ be a finite
  nilpotent group such that $\cc_G(x) = \{ [x,w] \mid w \in G \}$ for
  every $x \in G$.  Then $G$ has the `weak' Magnus
  property~\eqref{equ:magnus-f} for $f \colon \N \to \N$,
  $n \mapsto n$.
\end{lemma}

\begin{proof}
  Let $g,h \in G$.  Suppose that $k,l \in \N_0$ and
    \begin{equation} \label{equ:g-h-expressions}
    h = \prod\nolimits_{i=1}^k (g^{e_i})^{v_i} \qquad
    \text{and} \qquad g = \prod\nolimits_{j=1}^l (h^{d_j})^{w_j}
  \end{equation}
  for suitable $e_1, \ldots, e_k, d_1, \ldots, d_l \in \{1,-1\}$ and
  $v_1, \ldots, v_k, w_1, \ldots, w_l \in G$.  In particular, this
  implies that $\langle g \rangle^G = \langle h \rangle^G$.  By
  symmetry, it suffices to show that there exist $w \in G$ and
  $s \in [-l,l]_\Z$ such that $h = (g^s)^w$.

  If $g = 1$ then $h = 1$, and no further explanations are necessary.
  Now suppose that $g \ne 1$, and write
  $M = \cc_G(g) = \cc_G(h) \nt G$; see \cref{lem:cocentraliser}.  From
  \cref{equ:g-h-expressions} we deduce that $h \equiv_M g^s$, for some
  $s \in [-k,k]_\Z$.

  We claim that $\langle g \rangle = \langle g^s \rangle$.  For this
  it is enough to show that $p \nmid s$ for every prime $p$ that divides
  the order of~$g$.  The finite nilpotent group $G$ is the direct
  product of its Sylow subgroups; let $\bar x$ denote the image of
  $x \in G$ under the canonical projection onto the Sylow $p$-subgroup
  of~$G$.  Then $\bar g \ne 1$ implies that
  $g \in \mathrm{Z}_i(\overline{G}) \smallsetminus
  \mathrm{Z}_{i-1}(\overline{G})$, for suitable $i \in \N$, hence
  $\overline{M} \le \mathrm{Z}_{i-1}(\overline{G})$ and
  $\bar g \not\in \overline{M}$.  Consequently,
  $\langle \bar g \rangle \overline{M} = \langle \bar g
  \rangle^{\overline{G}} = \langle \bar h \rangle^{\overline{G}} =
  \langle {\bar g}^s \rangle \overline{M}$ implies $p \nmid s$.

  Using our general assumption on cocentralisers in~$G$, we deduce
  from $\langle g \rangle = \langle g^s \rangle$ that
  $M = \cc_G(g^s) = \{ [g^s,w] \mid w \in G \}$.  Therefore
  $h \equiv_M g^s$ shows that there exists $w \in G$ such that
  $h = (g^s)^w$.
\end{proof}

It remains to verify that the condition on cocentralisers in
\cref{lem:commutator-set-closed} applies to the concrete groups
$G = G_p$, defined in~\eqref{equ:concrete-Gp}.  Recall that $p > 2$.
Actually, it is convenient to check the required property for the
compact $p$-adic analytic group
\[
H = \langle t, a \mid [a,t] = a^p \rangle_{\text{pro-$p$}} \cong
(1+p\Z_p) \ltimes \Z_p,
\]
where $\Z_p$ denotes the ring of $p$-adic integers and the
multiplicative group of one-units
$1+p\Z_p = \overline{\langle 1+p \rangle} = \{ (1+p)^\lambda \mid
\lambda \in \Z_p \}$ acts naturally on the additive group~$\Z_p$.  The
group $H$ maps naturally onto~$G$, with kernel
$\overline{\langle t^{p^{c-1}}, a^{p^c} \rangle}$, and it is easy to
see that the condition on cocentralisers that we are interested in is
inherited by factor groups.

Let $h \in H$.  For $h=1$, it is clear that
$\{ [h,y] \mid y \in H \} = \{1\}$ is closed under multiplication.
Now suppose that $h \ne 1$.  Then $h$ is of the form
$h = t^\lambda a^\mu$ for uniquely determined
$\lambda, \mu \in \Z_p$, not both equal to~$0$.  Easy
computations show:
\begin{align*}
  \{ [h,a^\nu] \mid \nu \in \Z_p \} %
  & = \{ a^{(-(1+p)^\lambda + 1)\nu} \mid \nu \in
    \Z_p \} = \{ a^\sigma \mid \sigma \in p^{1+v_p(\lambda)}
    \Z_p \}, \\
  \{ [h,t^\nu] \mid \nu \in \Z_p \} %
  & = \{ a^{(-1 + (1+p)^\nu) \mu} \mid \nu \in
    \Z_p \} = \{ a^\sigma \mid \sigma \in p^{1+v_p(\mu)}
    \Z_p \}, 
\end{align*}
where $v_p \colon \Z_p \to \N_0 \cup \{\infty\}$
denotes the $p$-adic valuation map.

Put $m = 1+\min \{ v_p(\lambda), v_p(\mu) \}$.  As
$A_m = \{ a^\sigma \mid \sigma \in p^m\Z_p \}$ is a closed
normal subgroup of~$H$, the well-known commutator identity
$[a,bc] = [a,c] [a,b]^c$ for arbitrary group elements $a,b,c$ shows
that
\[
  \{ [h,y] \mid y \in H \} = A_m \trianglelefteq_\mathrm{c} H
\]
is indeed closed under multiplication.

%%%%%


\begin{bibdiv}
  \begin{biblist}
    % \bib{Bo05}{article}{
    % author={Bogopolski, Oleg},
    % title={A surface groups analogue of a theorem of Magnus},
    % conference={
    % title={Geometric methods in group theory},
    % },
    % book={
    % series={Contemp. Math.},
    % volume={372},
    % publisher={Amer. Math. Soc., Providence, RI},
    % },
    % date={2005},
    % pages={59--69},
    % review={\MR{2139677}},
    % doi={10.1090/conm/372/06874},
    % }
    \bib{BoSv08}{article}{ %
      author={Bogopolski, Oleg}, %
      author={Sviridov, Konstantin}, %
      title={A Magnus theorem for some one-relator groups}, %
      conference={ %
        title={The Zieschang Gedenkschrift}, }, %
      book={ %
        series={Geom. Topol. Monogr.}, %
        volume={14}, %
        publisher={Geom. Topol. Publ., Coventry}, %
      }, %
      date={2008}, %
      pages={63--73}, %
      review={\MR{2484697}}, %
      doi={10.2140/gtm.2008.14.63}, %
    } %

    \bib{Bo60}{article}{ %
      author={Bovdi, A. A.}, %
      title={Group rings of torsion-free groups}, %
      journal={Sibirsk. Mat. \v{Z}.}, %
      volume={1}, %
      date={1960}, %
      pages={555--558}, %
      issn={0037-4474}, %
      review={\MR{0130919}}, %
    } %
 
    \bib{ChKe90}{book}{ %
      author={Chang, C. C.}, %
      author={Keisler, H. J.}, %
      title={Model theory}, %
      series={Studies in Logic and the Foundations of Mathematics}, %
      volume={73}, %
      edition={3}, %
      publisher={North-Holland Publishing Co., Amsterdam}, %
      date={1990}, %
      pages={xvi+650}, %
      isbn={0-444-88054-2}, %
      review={\MR{1059055}}, %
    } %

    \bib{ChDo10}{article}{ %
      author={Chillag, David}, %
      author={Dolfi, Silvio}, %
      title={Semi-rational solvable groups}, %
      journal={J. Group Theory}, %
      volume={13}, %
      date={2010}, %
      number={4}, %
      pages={535--548}, %
      issn={1433-5883}, %
      review={\MR{2661654}}, %
      doi={10.1515/JGT.2010.004}, %
    } %
                 
    \bib{ClMaZy17}{book}{ %
      author={Clement, Anthony E.}, %
      author={Majewicz, Stephen}, %
      author={Zyman, Marcos}, %
      title={The theory of nilpotent groups}, %
      publisher={Birkh\"{a}user/Springer, Cham}, %
      date={2017}, %
      pages={xvii+307}, %
      isbn={978-3-319-66211-4}, %
      isbn={978-3-319-66213-8}, %
      review={\MR{3729243}}, %
    } %
                 
    % \bib{Ed89}{article}{
    % author={Edjvet, M.},
    % title={A Magnus theorem for free products of locally indicable
    % groups},
    % journal={Glasgow Math. J.},
    % volume={31},
    % date={1989},
    % number={3},
    % pages={383--387},
    % issn={0017-0895},
    % review={\MR{1021812}},
    % doi={10.1017/S001708950000793X},
    % }

    \bib{Fe19}{article}{ %
      author={Feldkamp, Carsten}, %
      title={A Magnus theorem for some amalgamated products}, %
      journal={Comm. Algebra}, %
      volume={47}, %
      date={2019}, %
      number={12}, %
      pages={5348--5360}, %
      issn={0092-7872}, %
      review={\MR{4019345}}, %
      doi={10.1080/00927872.2019.1623235}, %
    } %
                 
    \bib{Fe21}{article}{ %
      author={Feldkamp, Carsten}, %
      title={A Magnus extension for locally indicable groups}, %
      journal={J. Algebra}, %
      volume={581}, %
      date={2021}, %
      pages={122--172}, %
      issn={0021-8693}, %
      review={\MR{4249758}}, %
      doi={10.1016/j.jalgebra.2021.04.011}, %
    } %

    \bib{Gr57}{article}{ %
      author={Gruenberg, K. W.}, %
      title={Residual properties of infinite soluble groups}, %
      journal={Proc. London Math. Soc. (3)}, %
      volume={7}, %
      date={1957}, %
      pages={29--62}, %
      issn={0024-6115}, %
      review={\MR{87652}}, %
      doi={10.1112/plms/s3-7.1.29}, %
    } %
    
    \bib{HiNeNe49}{article}{ %
      author={Higman, Graham}, %
      author={Neumann, B. H.}, %
      author={Neumann, Hanna}, %
      title={Embedding theorems for groups}, %
      journal={J. London Math. Soc.}, %
      volume={24}, %
      date={1949}, %
      pages={247--254}, %
      issn={0024-6107}, %
      review={\MR{32641}}, %
      doi={10.1112/jlms/s1-24.4.247}, %
    } %
                 
    \bib{KlKu16}{article}{ %
      author={Klopsch, Benjamin}, %
      author={Kuckuck, Benno}, %
      title={The Magnus property for direct products}, %
      journal={Arch. Math. (Basel)}, %
      volume={107}, %
      date={2016}, %
      number={4}, %
      pages={379--388}, %
      issn={0003-889X}, %
      review={\MR{3552215}}, %
      doi={10.1007/s00013-016-0939-6}, %
    } %

    \bib{Ku82}{article}{%
      author={Kuz'min, Yu. V.},%
      title={Finite-order elements in free groups of some manifolds},%
      journal={Mat. Sb. (N.S.)},%
      volume={119(161)},%
      date={1982},%
      number={1},%
      pages={119--131, 160},%
      issn={0368-8666},%
      review={\MR{672413}},%
    }%

    \bib{Ma30}{article}{ %
      author={Magnus, Wilhelm}, %
      title={\"{U}ber diskontinuierliche Gruppen mit einer
        definierenden Relation. %
        (Der Freiheitssatz)}, %
      journal={J. Reine Angew. Math.}, %
      volume={163}, %
      date={1930}, %
      pages={141--165}, %
      issn={0075-4102}, %
      review={\MR{1581238}}, %
      doi={10.1515/crll.1930.163.141}, %
    } %

    \bib{MySo09}{article}{ %
      author={Myasnikov, A. G.},%
      author={Sokhrabi, M.},%
      title={Groups elementarily equivalent to a free 2-nilpotent
        group of finite rank},%
      journal={Algebra Logika},%
      volume={48},%
      date={2009},%
      number={2},%
      pages={203--244, 283--284, 286},%
      issn={0373-9252},%
      translation={%
        journal={Algebra Logic},%
        volume={48},%
        date={2009},%
        number={2},%
        pages={115--139},%
        issn={0002-5232},%
      },%
      review={\MR{2573019}},%
      doi={10.1007/s10469-009-9047-z},%
    }%
		
    \bib{MySo11}{article}{ %
      author={Myasnikov, Alexei G.},%
      author={Sohrabi, Mahmood},%
      title={Groups elementarily equivalent to a free nilpotent group
        of finite rank},%
      journal={Ann. Pure Appl. Logic},%
      volume={162},%
      date={2011},%
      number={11},%
      pages={916--933},%
      issn={0168-0072},%
      review={\MR{2817564}},%
      doi={10.1016/j.apal.2011.04.003},%
    }%
		
    
    \bib{Os10}{article}{ %
      author={Osin, Denis}, %
      title={Small cancellations over relatively hyperbolic groups and
        embedding theorems}, %
      journal={Ann. of Math. (2)}, %
      volume={172}, %
      date={2010}, %
      number={1}, %
      pages={1--39}, %
      issn={0003-486X}, %
      review={\MR{2680416}}, %
      doi={10.4007/annals.2010.172.1}, %
    } %
 
    \bib{Pa77}{book}{ %
      author={Passman, Donald S.}, %
      title={The algebraic structure of group rings}, %
      series={Pure and Applied Mathematics}, %
      publisher={Wiley-Interscience [John Wiley \& Sons], New
        York-London-Sydney}, %
      date={1977}, %
      pages={xiv+720}, %
      isbn={0-471-02272-1}, %
      review={\MR{0470211}}, %
    } %

    \bib{St89}{article}{%
      author={St\"{o}hr, Ralph},%
      title={On elements of order four in certain free central
        extensions of groups},%
      journal={Math. Proc. Cambridge Philos. Soc.},%
      volume={106},%
      date={1989},%
      number={1},%
      pages={13--28},%
      issn={0305-0041},%
      review={\MR{994077}},%
      doi={10.1017/S0305004100067955},%
    }%
 
    \bib{Wi10}{article}{ %
      AUTHOR = {Wilson, John S.}, %
      TITLE = {Free subgroups in groups with few relators}, %
      journal = {Enseign. Math. (2)}, %
      VOLUME = {56}, %
      YEAR = {2010}, %
      NUMBER = {1-2}, %
      PAGES = {173--185}, %
      ISSN = {0013-8584}, %
      review = {\MR{2674858}}, DOI = {10.4171/LEM/56-1-6}, %
    }
  \end{biblist}
\end{bibdiv}
\end{document}